\documentclass[aos,MSNbibl,nameyear,dvips]{arximspdf}
\usepackage{subenv}
\usepackage{graphicx}

\doi{10.1214/11-AOS889}
\volume{39}
\issue{4}
\pubyear{2011}
\firstpage{1916}
\lastpage{1962}

\makeatletter

\def\bptnote#1{}

\def\@bmisc[#1]{%
  \get@battribute{unstr}%
  \common@pub@types%
  \let\bauthor\bbl@bauthor%
  \let\bhowpublished\@firstofone%
  \def\borganization##1{{\bauthor@style ##1}}%
}

\newtheorem{prop}{Proposition}
\newtheorem{theorem}{Theorem}
\newtheorem{cor}{Corollary}
\newtheorem{lem}{Lemma}

\newcommand{\iint}{\int\!\!\int}
\newcommand{\iiint}{\int\!\!\int\!\!\int}

\newcommand{\DP}{\mathsf{DP}}

\newcommand{\Np}{{\mathcal{N}}}
\newcommand{\Nt}{{\tilde{\mathcal{N}}}}
\newcommand{\Eone}{{\mathrm{E}_1}}
\newcommand{\RA}{{\bbR\times A}}
\newcommand{\RK}{{\bbR\times K}}
\newcommand{\RO}{{\bbR\times\Omega}}
\newcommand{\gG}{{\phi}_{G}}
\newcommand{\gL}{{\phi}_{L}}
\newcommand{\gE}{{\phi}_{E}}
\newcommand{\gH}{{\phi}_{H}}
\newcommand{\Lmea}{\mathcal{L}}
\newcommand{\loc}{\chi}
\newcommand{\G}{{g}}

\newcommand{\GB}{{g}}
\newcommand{\Cov}{{\mathsf{Cov}}}
\newcommand{\E}{{\mathsf{E}}}
\newcommand{\V}{{\mathsf{Var}}}
\newcommand{\bY}{\mathbf{Y}}

\newcommand{\bbet}{\bolds{\beta}}
\newcommand{\bw}{\bolds{\omega}}
\newcommand{\dbdo}{d\beta\,d\omega}

\newcommand{\gd}{\dot\gamma}
\newcommand{\RRT}{{\bbR\times[1, \infty)\times\bbT}}

\newcommand{\bbW}{{\mathbb{W}}}
\newcommand{\bbB}{{\mathbb{B}}}
\newcommand{\Sob}{{\bbW^s_2}}
\newcommand{\Besov}{{\bbB^s_{pq}}}
\newcommand{\Besper}{{\bbB^{s *}_{pq}}}

\newcommand{\sY}{_{Y}}
\newcommand{\pat}{{\frac{\pi\alpha}{2}}}
\newcommand{\tpi}{{\frac{2}{\pi}}}

\newcommand{\km}{\mbox{ km}}

\newcommand{\eps}{\varepsilon}

\newcommand{\iI}{_{i\in I}}
\newcommand{\jJ}{_{0\le j< J}}
\newcommand{\jJe}{_{0\le j< J_\eps}}

\newcommand{\ind}{\stackrel{\mathrm{ind}}{\sim}}
\newcommand{\iid}{\stackrel{\mathrm{i.i.d.}}{\sim}}

\newcommand{\cmp}{h(\beta)}
\newcommand{\cmpj}{h(\beta_j)}
\newcommand{\stdcmp}{h_0(\beta)}
\newcommand{\half}{{\frac12}}
\newcommand{\thalf}{{\frac12}}

\newcommand{\bbN}{{\mathbb{N}}}
\newcommand{\bbR}{{\mathbb{R}}}
\newcommand{\cG}{\mathcal{G}}
\newcommand{\cH}{{\mathcal{H}}}
\newcommand{\cS}{\mathcal{S}}
\newcommand{\cT}{\mathcal{T}}
\newcommand{\cX}{\mathcal{X}}
\newcommand{\cZ}{\mathcal{Z}}

\newcommand{\Ga}{{\mathsf{Ga}}}

\newcommand{\Lv}{\mbox{{\textsf{L\'{e}vy}}}}
\newcommand{\NB}{{\mathsf{NB}}}
\newcommand{\No}{{\mathsf{No}}}
\newcommand{\Po}{{\mathsf{Po}}}
\newcommand{\St}{{\mathsf{St}}}
\newcommand{\Ca}{{\mathsf{C}}}
\newcommand{\Un}{{\mathsf{Un}}}

\newcommand{\ddd}{d}
\newcommand{\Scale}{\Lambda}
\newcommand{\LamS}{{\cS^d_+}}
\newcommand{\RS}{{\bbR\times\LamS}}
\newcommand{\Rd}{{\bbR^d}}

\newcommand{\nubw}{{\tilde\nu}} 

\newcommand{\uu}{{(\lambda(x - \loc))}}
\newcommand{\uuj}{{(\lambda_j(x-\loc_j))}}
\newcommand{\bbT}{{\mathbb{T}}}
\newcommand{\ui}{{[0,1)}}

\newcommand{\Z}{z} 
\newcommand{\Rp}{\Scale} 
\newcommand{\Ru}{\Scale} 

\makeatother

\begin{document}
\begin{frontmatter}

\title{Stochastic expansions using continuous dictionaries: L{\'{e}}vy adaptive regression kernels}
\runtitle{L{\'{e}}vy adaptive regression kernels}

\begin{aug}
\author[A]{\fnms{Robert L.} \snm{Wolpert}\corref{}\thanksref{t1}\ead[label=e1]{wolpert@stat.duke.edu}},
\author[A]{\fnms{Merlise A.} \snm{Clyde}\thanksref{t2}\ead[label=e2]{clyde@stat.duke.edu}%
\ead[label=u1,url]{http://www.stat.duke.edu}}
and
\author[B]{\fnms{Chong} \snm{Tu}\ead[label=e3]{Chong.Tu@pimco.com}}

\runauthor{R. L. Wolpert, M. A. Clyde and C. Tu}
\affiliation{Duke University, Duke University and PIMCO}
\address[A]{R. L. Wolpert\\
M. A. Clyde\\
Department of Statistical Science\\
Duke University\\
Durham, North Carolina 27708-0251\\
USA\\
\printead{e1}\\
\hphantom{E-mail: }\printead*{e2}\\
\printead{u1}}
\address[B]{C. Tu\\
PIMCO\\
1345 Avenue of the Americas\\
New York, New York 10105-4800\\
USA\\
\printead{e3}}
\end{aug}

\thankstext{t1}{Supported by NSF Grants
DMS-07-57549, PHY-09-41373 and NASA AISR Grant NNX09AK60G.}

\thankstext{t2}{Supported by NSF Grants DMS-0342172 and DMS-04-06115.}

\received{\smonth{5} \syear{2010}}
\revised{\smonth{1} \syear{2011}}

%
\begin{abstract}
This article describes a new class of prior distributions for nonparametric function estimation.
The unknown function is modeled as a limit of weighted sums of kernels or generator functions
indexed by continuous parameters that control local and global features such as their translation,
dilation, modulation and shape. L{\'{e}}vy random fields and their stochastic integrals are
employed to induce prior distributions for the unknown functions or, equivalently, for the number
of kernels and for the parameters governing their features. Scaling, shape, and other features of
the generating functions are location-specific to allow quite different function properties in
different parts of the space, as with wavelet bases and other methods employing overcomplete
dictionaries. We provide conditions under which the stochastic expansions converge in specified
Besov or Sobolev norms. Under a Gaussian error model, this may be viewed as a sparse regression
problem, with regularization induced via the L{\'{e}}vy random field prior distribution. Posterior
inference for the unknown functions is based on a reversible jump Markov chain Monte Carlo algorithm.
We compare the L{\'{e}}vy Adaptive Regression Kernel (LARK) method to wavelet-based methods using
some of the standard test functions, and illustrate its flexibility and adaptability in nonstationary
applications.
\end{abstract}

%


\end{frontmatter}

\section{Introduction}\label{sintr}
Popular approaches for nonparametric Bayesian estimation of unobserved
functions generally employ as prior distributions either Gaussian processes
(or random fields, in two or more dimensions) or mixtures of Dirichlet
processes. In this article, we focus attention on a wider class of processes,
L{\'{e}}vy random fields and their stochastic integrals.
These include Gaussian random fields as a limiting case, while Dirichlet
processes may be represented as ``normalized'' variants of the Gamma
L{\'{e}}vy
random field; L{\'{e}}vy random fields thus provide an important link
between two
of the random processes that form the foundation of Bayesian nonparametric
methods (see Section \ref{sother}).
In this article, we construct prior distributions for the mean function in
nonparametric regression as stochastic integrals of L{\'{e}}vy random fields.
Under suitable regularity, these can be expressed as stochastic expansions
using continuous dictionaries, permitting tractable Bayesian inference.
While our focus is on nonparametric regression, we hope that the reader will
see the possibilities of using L{\'{e}}vy random fields in other
contexts.\looseness=-1

To begin, suppose we have noisy measurements $\{Y_i\}\iI$ of an unknown
real-valued function $f\dvtx\cX\to\bbR$ observed at points $\{x_i\}\iI$
in some
complete separable metric space $\cX$, with $\E[Y_i] = f(x_i)$. In
nonparametric regression models, the mean function $f(\cdot)$ is often
regarded as an element of some Hilbert space $\cH$ of real-valued functions
on $\cX$, and is expressed as a linear combination of basis functions
$\{g_j\} \subset\cH$:
%
%
\begin{equation}\label{eexpan1}
f(x_i) = \sum\jJ g_j(x_i) \beta_j
\end{equation}
with some (finite or infinite) number $J$ of unknown coefficients $\{
\beta_j
\}\jJ$. There is a vast literature on classical and Bayesian
approaches for
estimating $f$ from noisy data using such methods as regression splines,
Fourier expansions, wavelet expansions, and kernel methods, including kernel
regression and support (or relevance) vector machines [see \citet
{ChuMarr1991},
\citet{CrisShaw2000}, \citet{DeniHolmMallSmit2002},
\citet{Vida1999}, \citet{Wahb1992}, for background and
references]. Many approaches,
including smoothing splines and support vector machines, use as many basis
elements, $J$, as there are data points, $n=|I|$, but employ regularization
to avoid over-fitting. Sparser solutions (using fewer basis elements,
$J\ll
n$) may be obtained through more stringent regularization penalties, as in
the Lasso [\citet{Tibs1996}] and Dantzig Selector [\citet
{CandTao2007}]
approaches, or (often equivalently) in Bayesian methods through choice of
prior distributions, as in relevance vector machines [\citet{Tipp2001}].
Sparse solutions may also be achieved by using variable selection techniques
to choose a few well-placed basis functions, perhaps in conjunction with
regularization [\citet{ChenDonoSaun1998}, \citet{DeniMallSmit1998b},
\citet{DiMaGenoKass2001}, \citet{MallZhan1993},
\citet{JohnSilv2005a}, \citet{SmitKohn1996},
\citet{WolfGodsNg2004}].

In most signal processing and other nonstationary applications, no single
(especially orthonormal) basis will lead to a sparse representation
[\citet{DonoElad2003}, \citet{WolfGodsNg2004}]. Overcomplete
dictionaries and frames\vadjust{\eject}
[\citet{Daub1992}, \citet{MallZhan1993}] provide larger
collections of
generating
elements $\{g_\omega\} _{\omega\in\Omega}$ than would a single
basis for
$\cH$, potentially allowing for more effective signal extraction and data
compression. Examples of overcomplete dictionaries include unions of bases,
Gabor frames, nondecimated or translational invariant wavelets, wavelet
packets, or more general kernel functions or generating functions $g(x,
\omega)$ where $\omega\in\Omega$ controls features (local or
global) of the
generating function, such as translations, dilations, modulations and shapes.
Because of the redundancy inherent in overcomplete representations,
coefficients for expansions using overcomplete dictionaries are not uniquely
determined. This lack of uniqueness is advantageous, permitting more
parsimonious representations from the dictionary than those obtained using
any single basis.\vadjust{\goodbreak}

In this article, we develop a fully Bayesian method for the sparse regression
problem using stochastic expansions [\citet{AbraSapaSilv2000}] of
continuous dictionaries. We begin in Section \ref{spriors} by introducing
L{\'{e}}vy
random fields, which are used to induce prior distributions for $f \in
\cH$
through stochastic integration of a kernel function with respect to a~signed
infinitely divisible random measure. We call the new model class
\textit{L{\'{e}}vy Adaptive Regression Kernel} or ``LARK'' models. The LARK
framework allows both the number of kernels and kernel-specific
parameters to
adapt to any nonstationary features of $f$. Both finite and infinite
expansions are considered. Exploiting the construction of L{\'{e}}vy random
fields through Poisson random fields, we develop finite approximations to
infinite expansions in Section \ref{sapprox} that permit tractable
inference. In Section~\ref{sspaces}, we provide conditions under which the
functions are almost
surely in the same function space as the generating kernel. We describe the
hierarchical representations of LARK models in Section \ref{sspreg}
that enable
posterior inference for the LARK model using reversible jump Markov chain
Monte Carlo (RJ-MCMC) methods. In Section \ref{sother}, we discuss
relationships
among LARK and other popular parametric and nonparametric methods. We then
compare our LARK method to other procedures using simulated data in
Section~\ref{ssim} and real data in Section~\ref{sapps}. In Section
\ref{sdisc}, we
discuss possible extensions of the LARK model.

\section{Stochastic expansions and prior distributions} \label{spriors}

$!\!\!\!\!$To make inference about the unknown mean function $f\in\cH$
given noisy observations $Y_i$ of $f(x_i)$ for \mbox{$\{x_i\}\subset\cX$},
we must
first propose a prior distribution on $\cH$ for $f$.
Let $\Omega$ be a complete separable metric space and
${\phi}\dvtx\cX\times\Omega\to\bbR$ a Borel measurable function, and
set ${\phi}_j(x_i)
\equiv{\phi}(x_i, \omega_j)$ for some collection
$\{\omega_j\}\subset\Omega$.
As a slight extension of the basis expansion of (\ref{eexpan1}), set
%
%
\begin{equation}\label{eexpan2}
f(x) \equiv\sum\jJ{\phi}(x, \omega_j) \beta_j
\end{equation}
for a \textit{random} number $J\le\infty$ of randomly drawn pairs
$(\beta_j,\omega_j)\in\bbR\times\Omega$. This is equivalent to
specifying a
random signed Borel measure $\Lmea(d \omega) =\sum\beta_j\delta
_{\omega_j}
(d\omega)$ on $\Omega$, giving the equivalent representation:
%
%
\begin{equation}\label{efint}
f(x) = \int_{\Omega} {\phi}(x, \omega) \Lmea(d \omega).
\end{equation}
The task of assigning prior distributions to functions $f(\cdot)$ of
the form
(\ref{eexpan2}) is equivalent to that of specifying prior
distributions for
the random measure~$\Lmea(d \omega)$ in (\ref{efint}), that is, to
specifying
consistent joint probability distributions for all random vectors of
the form
$(\Lmea(A_1), \ldots, \Lmea(A_k))$ for disjoint Borel sets
$A_i\subset\Omega$. L{\'{e}}vy random measures, those for which $\{
\Lmea
(A_i)\}$
are independent for disjoint $\{A_i\}$, are ideal for this purpose,
since (as
we will see in Section \ref{sspost}) they are simple to construct and
amenable to
posterior simulation. To make ideas more concrete, we first describe
possible choices for the generating functions ${\phi}(x, \omega)$
used in our
stochastic expansions and then proceed with the presentation of L{\'
{e}}vy random
measures in Section~\ref{sslrm}.

\subsection{Generating functions}\label{sskern}

Possible choices for ${\phi}(x, \omega)$ for $\cX=\bbR$ include
translation-invariant kernel functions, such as the Gaussian
%
%
\begin{subequation}\label{ekernels}
\begin{eqnarray}\label{eGaussian}
\gG(x, \omega) &\equiv&
\exp\bigl\{-{\tfrac12} \lambda(x - \loc)^2 \bigr\}
\end{eqnarray}
or the Laplace
\begin{eqnarray}\label{eLaplace}
\gL(x, \omega) &\equiv&
\exp\{- \lambda|x - \loc|\}
\end{eqnarray}
kernels with $\omega\equiv(\loc, \lambda) \in\cX\times\bbR^+
\equiv\Omega$. There is no need to restrict attention to symmetric
(e.g., Mercer) kernels, as required in the conventional Support Vector
Machine (SVM) approach [\citet{LawKwok2001},
\citet{Soll2002}]. Asymmetric kernels, such as the one-sided
exponential
\begin{eqnarray}\label{eExponential}
\gE(x,\omega) &\equiv&
\exp\{-\lambda(x-\loc)\} \mathbf{1}_{\{x>\loc\}}
\end{eqnarray}
are useful, for example, in modeling pollutant dissipation over time.
Other possibilities include piecewise-constant Haar wavelets on
$\cX=(0,1]$,
\begin{eqnarray}\label{eHaar}
\gH(x, \omega) &\equiv&
\mathbf{1}_{\{0 < \lambda(x - \loc) \le1\}}
\end{eqnarray}
or continuous rescaling and shifting of other wavelet functions
\begin{eqnarray}\label{eWave}
{\phi}_{\psi}(x, \omega) &\equiv&
\lambda^{1/2}\psi\bigl(\lambda(x - \loc)\bigr).
\end{eqnarray}
In each of these examples, $\Omega$ is a location-scale space with
location parameter~$\loc$ and parameter $\lambda$ determining the scale.
Higher-dimensional spaces~$\cX$ may be accommodated in a similar way; for
example, in Section \ref{ssair} we use space--time kernel
\begin{eqnarray}\label{eST}
{\phi}_{\mathrm{ST}}(x,\omega)
&\equiv&
\exp\bigl\{-\tfrac12(s-\sigma)'\Lambda(s-\sigma)-\lambda|t-\tau
|\bigr\}
\end{eqnarray}
for space--time point $x=(s,t)\in\bbR^2\times\bbR_+$;
here $\omega=(\sigma,\tau,\Lambda,\lambda)$ includes a~spa\-ce--time point
$(\sigma,\tau)\in\bbR^2\times\bbR_+$, a positive-definite spatial
dispersion matrix $\Lambda\in\cS^+_2$, and a temporal decay rate
$\lambda\in\bbR_+$.
\end{subequation}

\subsection{\texorpdfstring{L{\'{e}}vy random measures}{Levy random measures}}\label{sslrm}

$\!\!\!\!\!$For any $\nu^+\ge0$ and any probability distribu\-tion $\pi(\dbdo)$ on
$\RO$, let $J\,{\sim}\,\Po(\nu^+)$ be Poisson-distributed with~mean~$\nu^+$, and let
$\{(\beta_j,\omega_j)\}\jJ\iid\pi(\dbdo)$; then the random measure
given by
%
%
\begin{equation}\label{epmass}
\Lmea(A)\equiv\sum\jJ\mathbf{1}_{A}(\omega_j) \beta_j
\end{equation}
assigns independent infinitely-divisible (henceforth ``ID'') random
variables $\Lmea(A_i)$ 
to disjoint Borel sets $A_i\subset\Omega$, with characteristic functions
%
%
\begin{equation}\label{erf}
\E\bigl[e^{{i t \Lmea(A)}}\bigr] =
\exp\biggl\{\iint_{\RA} (e^{i t \beta}-1)
\nu(\dbdo)\biggr\}
\end{equation}
with $\nu(\dbdo)\equiv\nu^+ \pi(\dbdo)$.
More generally, the ``L{\'{e}}vy measure'' $\nu(\dbdo)$ need not be
finite for the random measure $\Lmea$ to be well defined, so long as the
integral in (\ref{erf}) converges for all $t\in\bbR$; since the
integrand is
bounded on all of $\RO$ and is of order $O(\beta)$ near $\beta
\approx0$,
this will hold for any measure that satisfies the {local $L_1$}
integrability condition
%
%
\begin{equation}\label{el1bound}
\iint_{\RK} (1\wedge|\beta| )
\nu(\dbdo)<\infty
\end{equation}
for each compact $K\subset\Omega$. The mean and variance, when they exist,
are given by $\E[\Lmea(A)] = \iint_{\RA}\beta\nu(\dbdo)$ and
$\V[\Lmea(A)] = \iint_{\RA} \beta^{ 2} \nu(\dbdo)$,
respectively.

\citet{KhinLevy1936} showed that the most general ID random variables
[and hence the most general ID-valued random measures; see \citet{RajpRosi1989},
Proposition 2.1] have characteristic functions of the form
%
%
\begin{eqnarray}\label{elevy-full}
\E\bigl[e^{i t \Lmea(A)}\bigr] &=&
\exp\biggl\{i t \delta(A) - {\thalf} t^2 \Sigma(A)\nonumber\\[-8pt]\\[-8pt]
&&\hphantom{\exp\biggl\{}
{}+ \iint_{\RA} \bigl(e^{i t \beta}-1
- i t \stdcmp\bigr) \nu(\dbdo)\biggr\},\nonumber
\end{eqnarray}
where $\stdcmp\equiv\beta\mathbf{1}_{[-1,1]}(\beta)$, determined
uniquely by
the characteristic triplet of sigma-finite measures $(\delta, \Sigma,
\nu)$
consisting of a signed measure $\delta(d \omega)$ and a positive measure
$\Sigma(d\omega)$ on $\Omega$, and a positive measure
$\nu(\dbdo)$ on $\bbR\times\Omega$ that satisfies the {local
$L_2$} integrability condition
%
%
\begin{equation}\label{el2bound}
\iint_\RK(1 \wedge\beta^2) \nu(\dbdo) < \infty
\end{equation}
for each compact $K\subset\Omega$ and $\nu(\{0\}, \Omega) = 0$ (for more
details on this nonstationary version of the classic L{\'{e}}vy--Khinchine
formula see Jacod and Shiryaev [(\citeyear{JacoShir1987}), page 75],
Cont and Tankov [(\citeyear{ContTank2004}), pages 457--459]
or \citet{WolpTaqq2005}).

The role of the \textit{compensator} function $\stdcmp$ is to make
the last integrand in (\ref{elevy-full}) bounded and $O(\beta^2)$ near
$\beta\approx0$, permitting the replacement of (\ref{el1bound})
with the
weaker condition (\ref{el2bound}); in this case $\Lmea(d\omega)$
may have
countably-many points of support $\{\omega_j\}\subset\Omega$ whose
magnitudes $\{\beta_j\}$ are not absolutely summable, precluding a
representation of the form (\ref{epmass}). The compensator $\stdcmp$
may be
replaced by any bounded measurable function satisfying
%
%
\begin{equation}\label{ecomp}
\cmp= \beta+ O(\beta^2), \qquad\beta\approx0,
\end{equation}
with a corresponding replacement of $\delta(d\omega)$ with $\delta_{h}
(d\omega)=\delta(d\omega) +\int_{\bbR} [\cmp-\stdcmp]
\nu(\dbdo)$. Whenever (\ref{el1bound}) is satisfied, we may take
$\cmp\equiv0$ with the same adjustment to $\delta_0$.


By (\ref{elevy-full}) the random measure $\Lmea$ may be written as the
sum of
two independent parts: a Gaussian portion, assigning independent
normally-distributed random variables with mean $\delta_{h}(A_i)$ and
variance $\Sigma(A_i)$ to disjoint sets $A_i$, and the remaining portion,
with characteristic function
%
%
\begin{equation}\label{elevy-jump}
\E\bigl[e^{i t \Lmea(A)}\bigr] = \exp\biggl\{
\iint_{\RA} \bigl(e^{i t \beta} - 1
- i t \cmp\bigr) \nu(\dbdo)\biggr\}.
\end{equation}
We call a random signed measure $\Lmea$ with no Gaussian component
[i.e., an
ID-valued measure with $\Sigma(\Omega) = \delta_{h}(\Omega)=0$, that
satisfies (\ref{elevy-jump})] a \textit{L{\'{e}}vy random measure}.
Nonnegative
L{\'{e}}vy random measures satisfying (\ref{el1bound}) were called ``completely
random measures'' by \citet{King1967}.

\subsection{\texorpdfstring{L{\'{e}}vy random fields}{Levy random fields}}\label{sslrf}

A L{\'{e}}vy random measure $\Lmea$ satisfying (\ref{elevy-jump}) induces
a linear
mapping $\phi\mapsto\Lmea[\phi]$ from functions $\phi\dvtx\Omega\to
\bbR$ to
random variables $\Lmea[\phi] \equiv\int_\Omega\phi(\omega)
\Lmea
(d\omega)$; such a mapping\vspace*{1pt} is called a \textit{random field}.
For simple functions $\phi(\omega) =\sum a_i \mathbf{1}_{A_i}(\omega
)$ with each
$\bar A_i\subset\Omega$ compact, we set $\Lmea[\phi]\equiv\sum a_i
\Lmea(A_i)$
and verify that
%
%
\begin{equation}\label{elevy1}
\E\bigl[e^{i t \Lmea[\phi]}\bigr] = \exp\biggl\{
\iint_{\bbR\times\Omega} \bigl(e^{i t \phi(\omega) \beta}-1
- i t \phi(\omega) \cmp\bigr) \nu(\dbdo)\biggr\}.
\end{equation}
It is straightforward to extend this by continuity in probability to (at
least) all bounded measurable compactly-supported $\phi\dvtx\Omega\to
\bbR$. We
now present a~general construction based on Poisson random fields, the
key to
our approach to tractable posterior Bayesian inference.

\subsubsection{Poisson construction  {I}: Uncompensated}\label{sssPoiConsUns}

When $\nu(\dbdo)$ satisfies (\ref{el1bound}) (i.e., $|\beta|$ is locally
$\nu$-integrable at zero) we may take $\cmp\equiv0$ in (\ref
{elevy1}) and
construct $\Lmea$ as follows. Begin with a Poisson random measure
$\Np(\dbdo)\sim\Po(\nu)$ on $(\RO)$ that assigns independent
Poisson-distributed random variables $\Np(C_i)\,{\sim}\,\Po(\nu(C_i))$
with means $\nu(C_i)$ to disjoint Borel sets $C_i\,{\subset}\,(\RO)$. For any
Borel set $A\subset\Omega$ with compact closure $\bar A$ and
bounded\vadjust{\eject}
measurable compactly-supported $\phi\dvtx\Omega\to\bbR$, set $J\equiv
\Np(\bbR
\times A)$ and
%
%
\begin{eqnarray}\label{epois-unc}
\Lmea(A) &\equiv& \iint_{\bbR\times A} \beta\Np( \dbdo)
= \sum\jJ\mathbf{1}_A(\omega_j) \beta_j,\nonumber\\
\Lmea[\phi] &\equiv& \iint_\RO\beta\phi(\omega) \Np( \dbdo)
= \sum\jJ\phi(\omega_j) \beta_j,
\end{eqnarray}
where $\{(\beta_j,\omega_j)\}$ is the (random) set of $J\le\infty$ support
points of $\Np(\dbdo)$. The integrals and sums in (\ref{elevy1}),
(\ref{epois-unc})
are well defined for all $\phi$ for which
\[
\iint_{[-1,1]\times\Omega}
|\beta\phi(\omega)| \nu(\dbdo)<\infty,
\]
which by (\ref{el1bound}) includes
all bounded measurable compactly-supported \mbox{functions}.

For any Borel sets $A\subset\Omega$ and $B\subset\bbR$, the Poisson measure
$\Np$ assigns to the set $B\times A\subset\RO$ the number $\Np
(B\times A)$ of
$\Lmea$'s support points $\omega_j\in A$ with mass of sizes $\beta
_j\in B$.
By (\ref{el1bound}) this is necessarily finite if $A$ has compact
closure and
$B$ is bounded away from zero, but if $\nu(\RO)=\infty$ then $\Lmea
$ will
have $J=\infty$ support points in $\Omega$ altogether with (almost surely)
absolutely summable magnitudes $\sum\jJ\{|\beta_j|\dvtx\omega
_j\in A\}
<\infty$.

\subsubsection{Poisson construction  {II}: Compensated}\label{sssPoiConsCom}

The situation is more delicate in case the L{\'{e}}vy measure does not satisfy
(\ref{el1bound}), but only the weaker bound in (\ref{el2bound})
(i.e., if
$\beta^2$ is locally $\nu$-integrable but $|\beta|$ is not). Begin again
with the Poisson measure $\Np\sim\Po(\nu)$ on $\RO$, and introduce the
\textit{compensated} or centered Poisson measure $\Nt(\dbdo) \equiv
\Np(\dbdo) - \nu(\dbdo)$ with mean zero [\citet{Sato1999},
page 38], inducing
an isometry from $L_2(\RO,\nu(\dbdo))$ to the square-integrable
zero-mean random variables. Following \citet{WolpTaqq2005}, set
%
%
\begin{eqnarray} \label{epois-comp}
\Lmea(A) &\equiv& \iint_{\bbR\times A}
[\beta-\cmp] \Np( \dbdo)
+ \iint_{\bbR\times A}
\cmp\Nt( \dbdo), \nonumber\\
\Lmea[\phi] &\equiv& \iint_\RO
[\beta-\cmp] \phi(\omega) \Np( \dbdo)\\
&&{} + \iint_\RO
\cmp\phi(\omega) \Nt( \dbdo)\nonumber
\end{eqnarray}
for any measurable $\phi$ for which (\ref{epois-comp}) converges. If
(\ref{el1bound}) holds, one may simplify~(\ref{epois-comp}) to
%
%
\begin{subequation}\label{enonsense}
\begin{eqnarray}
\label{enonsensea}
\Lmea[\phi]&=&\iint_\RO\beta\phi(\omega) \Np(\dbdo)
-\iint_\RO h(\beta)\phi(\omega) \nu(\dbdo)\\
\label{enonsenseb}
&=&\sum\jJ\phi(\omega_j) \beta_j
+\delta_h[\phi]
\end{eqnarray}
\end{subequation}
showing that the role of the compensator is to add an $h$-dependent ``drift'' (or \textit{offset},
in higher dimensions) term $\delta_h[\phi]=-\iint_\RO h(\beta) \phi(\omega) \nu(\dbdo)$
to~(\ref{epois-unc}). When~(\ref{el1bound}) fails, however, both the uncompensated sum and
$\delta_h[\phi]$ in (\ref{enonsense}) will be infinite, while the representation of
(\ref{epois-comp}) remains valid under the following
conditions.\vspace*{-3pt}
\begin{theorem}\label{texist}
Let $\nu$ be a L{\'{e}}vy measure on $\RO$ satisfying (\ref
{el2bound}). Then~$\Lmea[\phi]$ is well defined by (\ref{epois-comp}) with characteristic
function given by (\ref{elevy1}) for compensator $h_0(\beta)\equiv
\beta\mathbf{1}_{\{|\beta|\le1\}}$ if $\phi$ satisfies
%
%
\begin{subequation}\label{eexist}
\begin{eqnarray}\label{eexa}
\iint_{[-1,1]^c\times\Omega} \bigl(1\wedge|\beta\phi(\omega
)|\bigr)
\nu( \dbdo)&<&\infty,\\
\label{eexb}
\iint_{[-1,1]\times\Omega} \bigl(|\beta\phi(\omega)|\wedge|
\beta\phi(\omega)|^2\bigr)
\nu( \dbdo)&<&\infty.
\end{eqnarray}
If, in addition, $\phi$ satisfies
\begin{eqnarray}\label{eexc}
\iint_{\bbR\times\Omega} (1\wedge\beta^2) |\phi
(\omega)|
\nu( \dbdo)&<&\infty,
\end{eqnarray}
then $\Lmea[\phi]$ is well defined for any compensator $\cmp$ satisfying
(\ref{ecomp}).\vspace*{-3pt}
\end{subequation}
\end{theorem}
\begin{pf}
Under these conditions, the integrands of the compensated and uncompensated
Poisson integrals in (\ref{epois-comp}) are in the Musielak--Orlicz spaces
for which those integrals are well defined; see
Rajput and Rosi{\'{n}}ski [(\citeyear{RajpRosi1989}), page 9],
\citet{KwapWoyc1992}.\vspace*{-3pt}
\end{pf}

In particular:\vspace*{-3pt}
\begin{cor}
$\Lmea[\phi]$ is well defined with characteristic function (\ref{elevy1})
for any function $\phi$ satisfying
%
%
\begin{equation}\label{eeasy}
\iint_\RO(1\wedge\beta^2)
\bigl(|\phi(\omega)|\vee\phi^2(\omega)\bigr)
\nu(\dbdo)<\infty,
\end{equation}
including [by (\ref{el2bound})] all\vspace*{1pt} bounded measurable compactly-supported
$\phi$. Thus,\break $\Lmea(A)=\Lmea[\mathbf{1}_A]$ is always well defined
for any Borel
set $A\subset\Omega$ with compact closure $\bar A$.\vspace*{-3pt}
\end{cor}

Similarly:\vspace*{-3pt}
\begin{prop}
$\!\!\!$For a L{\'{e}}vy measure $\nu$ satisfying (\ref{el1bound}), take
$\cmp\equiv0$;~then
\[
\mbox{(\ref{epois-unc})}\hspace*{57.5pt} \Lmea[\phi] \equiv\iint_\RO\beta
\phi(\omega) \Np( \dbdo)
= \sum\jJ\phi(\omega_j) \beta_j\hspace*{57.5pt}
\]
[with $J\equiv\Np(\bbR\times\Omega)\le\infty$] is well defined
with characteristic
function (\ref{elevy1}) for any $\phi$ satisfying
%
%
\begin{equation}\label{enocomp}
\iint_{\bbR\times\Omega} \bigl(1\wedge|\beta\phi(\omega)|
\bigr)
\nu( \dbdo)<\infty.
\end{equation}
\end{prop}


\subsection{\texorpdfstring{Constructing L{\'{e}}vy kernel
integrals}{Constructing Levy kernel integrals}}

Denote by $\Phi$ the linear space of functions $\phi\dvtx\Omega\to\bbR
$ for which
$\Lmea[\phi]$ has been defined; we have seen that this includes at
least all
bounded measurable compactly-supported functions $\phi$. Denote by
$\cG$ the
linear space of measurable functions $g\dvtx\cX\to\Phi$, and simplify notation
by writing ``$g(x,\omega)$'' for $g(x)(\omega)$. Each of the generating
functions introduced in (\ref{ekernels}) lies in~$\cG$. For any
$g\in
\cG$, we
can construct a random function $f\dvtx\cX\to\bbR$ by
%
%
\begin{eqnarray}\label{ef-def}
f(x) &\equiv&\Lmea[g(x)]\\
&=& \iint_\RO g(x,\omega) [\beta-\cmp] \Np(\dbdo)\nonumber\\
&&{}+\iint_\RO g(x,\omega) \cmp\Nt(\dbdo)\nonumber\\
&=& \sum\jJ g(x,\omega_j) [\beta_j-\cmpj]\nonumber\\
&&{}+\iint_\RO g(x,\omega) \cmp\Nt(\dbdo)\quad
\mbox{or}\nonumber\\
\label{ef-nocomp}
&=& \sum\jJ g(x,\omega_j) \beta_j \qquad\mbox{if
(\ref{el1bound}) holds so compensation is unneeded.}\hspace*{-27pt}
\end{eqnarray}
Integer moments of $f(x)$ are easy to compute, when they exist, from the
characteristic function given in (\ref{elevy1}), for example:
%
%
\begin{subequation}\label{emoments}
\begin{eqnarray}
\label{emean}
\E\{ f(x)\} &=&
\iint_\RO{\phi}(x,\omega) [\beta-\cmp]
\nu(\dbdo),\\
\label{ecov}
\Cov\{f(x_1),f(x_2)\} &=&
\iint_\RO{\phi}(x_1,\omega) {\phi}(x_2,\omega) \beta^2
\nu(\dbdo).
\end{eqnarray}
\end{subequation}

\subsection{\texorpdfstring{Examples of L{\'{e}}vy measures}{Examples of Levy measures}}

We now consider some specific examples of L{\'{e}}vy random fields and the
corresponding kernel integrals. Familiar examples include Poisson, Gamma,
Cauchy and more generally $\alpha$-Stable random fields.

\subsubsection{Compound Poisson processes}\label{sssFin}
$\!\!\!\!$The simplest model to consider would be that of (\ref{eexpan2}), with
finite L{\'{e}}vy measure satisfying $\nu^+\equiv\nu(\bbR\times\Omega)<\infty$,
reproduced here:
\setcounter{equation}{1}
\begin{equation}
f(x) \equiv\sum\jJ{\phi}(x, \omega_j) \beta_j.
\end{equation}
\setcounter{equation}{21}
This\vspace*{-1pt} has a Poisson-distributed number $J\sim\Po(\nu^+)$ of terms whose
locations~$\omega_j$ and magnitudes $\beta_j$ are i.i.d. with an arbitrary
distribution $\{\beta_j, \omega_j\} \iid\pi(\dbdo)$, hence L\'
{e}vy measure of
the form $\nu(\dbdo) = \nu^+ \pi(\dbdo)$. The marginal
distribution of
$f(x)$ at each $x\in\cX$ is compound Poisson.\vspace*{-3pt}

\subsubsection{Gamma random fields}\label{sssGRF}
The L{\'{e}}vy measure for the Gamma random field is infinite but
satisfies the
strong local $L_1$ integrability condition (\ref{el1bound}), obviating
compensation; in the homogeneous case, it is
%
%
\begin{equation}\label{enu-gamma}
\nu(\dbdo)= \beta^{-1}
e^{-\beta\eta}\mathbf{1}_{\{\beta>0\}} \,d\beta\gamma(d\omega)
\end{equation}
for some $\sigma$-finite measure $\gamma(d\omega)$ on $\Omega$, giving
$\Lmea(A) \sim\Ga(\gamma(A), \eta)$ [with mean $\gamma(A)/\eta$]
for Borel
measurable $A \subset\Omega$ with $\gamma(A)<\infty$. Because $\nu
$ is
concentrated on $\bbR_+$, the mass $\beta_j$ at each of the Gamma random
measure's support points $\omega_j$ is positive, so all the
coefficients in
the expression $f(x)=\sum{\phi}(x,\omega_j) \beta_j$ are
nonnegative. With a
nonnegative generating function ${\phi}\in\cG$, this provides a
direct way to
construct nonnegative mean functions $f\ge0$ without having to
transform the
responses $\{Y_i\}$ as Gaussian methods would require. The mean $\E[f(x)]
=\eta^{-1}\int g(x,\omega)\gamma(d\omega)$ is available from (\ref
{emean}), as
is the covariance from~(\ref{ecov}).\vspace*{-3pt}

\subsubsection{Symmetric Gamma random fields}\label{sssSGRF}
A symmetric analogue of the Gamma random field (\ref{enu-gamma}) has
L{\'{e}}vy
measure
%
%
\begin{equation}\label{enu-gamma-symm}
\nu(\dbdo)= |\beta|^{-1}
e^{-|\beta| \eta} \,d\beta\gamma(d\omega)
\end{equation}
on all of $\RO$, leading to random variables $\Lmea(A)$ distributed
as the
difference of two independent $\Ga(\gamma(A), \eta)$
variables, with
characteristic function $\E[e^{it\Lmea(A)}]
=(1+t^2/\eta^2)^{-\gamma(A)}$. Both the standard positive Gamma random
measure and this symmetric version satisfy the local $L_1$ bound (\ref
{el1bound}), hence no compensation is required so we may take $\cmp
\equiv0$
and employ the simple construction (\ref{ef-nocomp}) of $f(x)$.
The mean $\E[f(x)]=0$ vanishes for the symmetric Gamma random field,
or for
any other L{\'{e}}vy random field with a symmetric (in $\pm\beta$)
L{\'{e}}vy
measure
satisfying (\ref{el1bound}). Covariances are available from (\ref{ecov}).
Nearly all of the commonly used isotropic geostatistical covariance functions
[see \citet{ChilDelf1999}, Section 2.5] may be achieved by the
choice of a
suitable generating kernel $g(x,\cdot)$ and L{\'{e}}vy measure
$\nu(\dbdo)$; see \citet{ClydWolp2007} for specific examples.
%

\vspace*{-3pt}\subsubsection{\texorpdfstring{Symmetric $\alpha$-Stable random fields}{Symmetric alpha-Stable random fields}}
\label{sssSaS}

Symmetric $\alpha$-Stable (\textsf{S$\alpha$S}) L{\'{e}}\-vy random
fields have L{\'{e}}vy measure
%
%
\begin{equation}\label{nu-sas}
\nu(\dbdo)=c_\alpha\alpha|\beta|^{-1 -\alpha} \,d\beta\,
\gamma(d\omega)
\end{equation}
%
on $\RO$ for some $0<\alpha<2$ and $\sigma$-finite positive measure
$\gamma(d\omega)$, where $c_\alpha= (1/\pi)\Gamma(\alpha)\sin
(\pi\alpha/2)$, giving $\Lmea(A)\sim\St(\alpha,0,\gamma(A),0)$
[in parametrization~(M) of \citet{Zolo1986}, page 11] with
infinite variance (and thus
no meaningful covariance function for $f(x)\equiv\Lmea[g(x)]$). This infinite
L{\'{e}}vy measure satisfies~(\ref{el2bound}) for all $0<\alpha<2$, but satisfies the
stronger local $L_1$ condition~(\ref{el1bound}) only for $0<\alpha<1$; thus
compensation is required to construct \textsf{S$\alpha$S} random fields with $1\le
\alpha<2$, including the\vadjust{\eject}
Cauchy case of $\alpha=1$. One can show that $f(x)$ is well defined for
any ${\phi}(x,\cdot)\in L_\alpha(\Omega,\gamma(d\omega))$, including
the generating functions of (\ref{ekernels}).
The \textsf{S$\alpha$S} fields have heavier tails than, for example,
the symmetric Gamma
fields of Section \ref{sssSGRF}, and may be more appropriate for problems
where one
might expect $f(\cdot)$ to include by a few heavily weighted kernels.

\vspace*{8pt}\section{Approximations for implementing kernel integrals}\label{sapprox}

Computer simulations of L{\'{e}}vy random measures $A\mapsto\Lmea(A)$
and random
fields $\phi\mapsto\Lmea[\phi]$ associated with \textit{finite}
L{\'{e}}vy
measures $\nu$ may be constructed as in (\ref{epmass}), (\ref{epois-unc}),
simply by
setting $\nu^+\equiv\nu(\RO)$ and drawing\vspace*{1pt} $J\sim\Po(\nu^+)$ and
$\{(\beta_j,\omega_j)\}\jJ\iid\pi(\dbdo) \equiv\nu(\dbdo)/\nu
^+$. If
$\nu(\RO)=\infty$ however the sums in these equations will include countably
infinitely-many terms, and may not be absolutely summable. We now construct
an approximating set of finite L{\'{e}}vy measures $\{\nu_\eps\}$
indexed by
$\eps>0$ and show that the approximate L{\'{e}}vy random fields
$\Lmea_\eps[\phi]$ converge to the random field $\Lmea[\phi]$
given in
(\ref{epois-comp}). Note that~$\eps$ is \textit{not} a model parameter.
It is
only a device used for two purposes: as a~tool in the theorems constructing
LARK models (in this section) and establishing their properties (in
Section \ref{sspaces}), and to enable the construction of practical numerical
methods to approximate LARK models within specified error bounds (in
Section \ref{sspreg}).
\begin{theorem}\label{tconverge}
Let $\nu$ be a L{\'{e}}vy measure defined on $\RO$ satisfying (\ref
{el2bound})
and $\phi\in\Phi$ satisfying (\ref{eexist}). Take $\{K_\eps\}$ to
be any
family of compact sets increasing to $\Omega$ as $\eps\to0$, and for any
Borel sets $A\subset\Omega$ and $B\subset\bbR$ and let $\nu_\eps$
be the
unique Borel measure on $\bbR\times\Omega$ satisfying
%
%
\begin{equation}\label{etrunc-nu}
\nu_\eps(B\times A) \equiv\nu\bigl((B\cap[-\eps,\eps]^c)
\times(A\cap K_\eps)\bigr)
\end{equation}
for $B\subset\bbR$, $A\subset\Omega$ [note $\nu_\eps^+\equiv
\nu_\eps(\bbR\times\Omega)<\infty$]. Let $h(\cdot)$ be any bounded
measurable compensator function on $\bbR$ satisfying $\cmp= \beta+
O(\beta^2)$ for $\beta$ near zero. Then as $\eps\to0$, the random
variables
%
%
\begin{eqnarray}\label{epoi-comp-eps}
\Lmea_\eps[\phi] &\equiv&
\iint_{[-\eps,\eps]^c \times K_\eps}
\beta\phi(\omega) \Np(\dbdo)\nonumber\\[-8pt]\\[-8pt]
&&{}- \iint_{[-\eps,\eps]^c \times K_\eps}
\cmp\phi(\omega) \nu(\dbdo)\nonumber
\end{eqnarray}
converge in probability to $\Lmea[\phi]$ of (\ref{epois-comp}).
\end{theorem}
\begin{pf}
The error in approximating $\Lmea[\phi]$ of (\ref{epois-comp}) by
$\Lmea_\eps[\phi]$ of (\ref{epoi-comp-eps}) is
%
%
\begin{eqnarray}\label{epoi-comp-res}
\Lmea[\phi]-\Lmea_\eps[\phi] &=&
\iint_{N_\eps}\bigl(\beta- \cmp\bigr) \phi(\omega)
\Np(\dbdo)\nonumber\\[-8pt]\\[-8pt]
&&{}+ \iint_{N_\eps} \cmp\phi(\omega) \Nt(\dbdo),\nonumber
\end{eqnarray}
where $N_\eps\equiv\{(\beta,\omega)\dvtx |\beta| \le\eps\mbox{ or
}\omega\in
K^c_\eps\}$. The first term in (\ref{epoi-comp-res}) converges to
zero almost
surely, and the second in $L_1$, as $\eps\to0$; see the \hyperref
[app]{Appendix}
for details.
\end{pf}

The approximation $\Lmea_\eps[\phi]$ is the sum of a L{\'{e}}vy random
field with
finite L{\'{e}}vy measure $\nu_\eps$ [hence with simple representation
(\ref{epois-unc})] and a deterministic drift term $\delta_\eps[\phi
]$ given by
the second integral in (\ref{epoi-comp-eps}). The drift vanishes whenever
$\nu(\dbdo)$ is symmetric in $\pm\beta$ and $\cmp$ is odd.

\begin{cor}\label{cfe}
If either
\textup{(a)} $\nu(\dbdo)$ satisfies (\ref{el1bound}), or
\textup{(b)} $\nu(\dbdo)$ satisfies (\ref{el2bound}) and is
even in $\pm\beta$, and also $\cmp$ is an odd function, then
for each $x \in\cX$,
%
%
\begin{equation}\label{eapprox-f}
f_\eps(x) \equiv\sum\jJe g(x,\omega_j) \beta_j
\end{equation}
with
\[
J_\eps\sim\Po(\nu^+_\eps),\qquad
\{\beta_j,\omega_j\}\jJe\mid J_\eps\stackrel{\mathit{i.i.d.}}{\sim}\nu_\eps
(\dbdo) /\nu^+_\eps
\]
converges to $f(x)$ in probability as $\eps\to0$.
\end{cor}
\begin{pf}
With $f_\eps(x) \equiv\Lmea_\eps[g(x)]$,
%
%
\begin{eqnarray}\label{egcomp}\quad
f_\eps(x) &=& \int_\Omega g(x,\omega)\Lmea_\eps(d\omega)
\nonumber\\[-8pt]\\[-8pt]
&=& \iint_{\bbR\times\Omega}
g(x,\omega)\beta\Np_\eps(\dbdo)
- \iint_{\bbR\times\Omega} g(x,\omega)
\cmp\nu_\eps(\dbdo)\nonumber
\end{eqnarray}
with $\Np_\eps(\dbdo) \sim\Po(\nu_\eps(\dbdo))$.
If $\nu$ satisfies (\ref{el1bound}), then without loss of generality
take the
compensator function $\cmp\equiv0$. In both cases (a) and (b), the second
integral in (\ref{egcomp}) vanishes, leading to (\ref{eapprox-f}) [cf.
(\ref{eexpan2})].
\end{pf}

Note that in case (b) 
the $\{g(x,\omega_j) \beta_j\}$ are not absolutely summable so
``$\sum_{j=0} ^\infty g(x,\omega_j) \beta_j$'' does not converge in the
Lebesgue sense. In each of our applications the conditions of Corollary
\ref{cfe} hold, allowing us to approximate~$\nu$ by a finite L{\'{e}}vy measure $\nu_\eps$
[and $\Lmea$ by $\Lmea_\eps\sim\Lv(\nu_\eps)$], and exploit the resulting
Poisson representation for inference.

%

\section{Function spaces for LARK models}\label{sspaces}

Theorem \ref{tconverge} and Corollary \ref{cfe} establish pointwise
convergence of $f_\eps(x)$ to $f(x)$ as $\eps\to0$; in this section we
provide conditions to ensure that $f_\eps(\cdot)\to f(\cdot)$ in appropriate
Besov or Sobolev norms if the generating functions lie in the same space.

For $s\ge0$ and 
$d\in\bbN$ denote by $\Sob(\bbR^d)$ the Sobolev space of real-valued
square-integrable functions $f(\cdot)\in L_2(\bbR^d)$
[\citet{Sobo1991}, Section 1.7, \citet{ReedSimo1975b},
page 50]
with finite
Sobolev norm
%
%
\begin{equation}\label{esob}
\|{f}\|_\Sob= \biggl\{\frac{1}{(2\pi)^d}
\int_{\bbR^d} (1+|\xi|^2)^s |\hat f(\xi)|^2
\,d\xi\biggr\}^{1/2}
\end{equation}
with Fourier transforms defined for $f\in L_1(\bbR^d)$ by
\[
\hat f(\xi) = \int_{\bbR^d} e^{i\xi\cdot x} f(x) \,dx
\]
and by $L_2$ limits for $f\in L_2(\bbR^d)$;
here $d\xi$ and $d x$ denote the Lebesgue volume element in $\bbR^d$,
and $\xi\cdot x$ denotes the Euclidean inner product.
Each $\Sob$ is a~Banach space, hence complete. By Plancherel's
theorem, each
$f\in\Sob$ with $s\ge0$ has $s$ distributional derivatives in
$L_2(\bbR)$,
and by Sobolev's lemma has~$k$ continuous derivatives for each integer
$0\le
k<s-d/2$.

Besov spaces constitute a flexible family that includes elements with wide
spatial irregularity. The Besov space\vspace*{1pt} $\Besov$ consists of those $f\in
L_p(\bbR^d)$ whose Besov semi-norms are finite. Several equivalent Besov
semi-norms appear in the literature [\citet{Trie1992},
Theorem 2.6.1, page 140]; we use the definition given as equation 2 of
that theorem. For
$p,q\ge0$ and $s> d(1/p - 1)_+$ and for any integer $m> s$
($m=1+\lfloor
s\rfloor$ is easiest), set
\[
|f|^s_{pq} = \biggl(\int_{|h|\le1} |h|^{-sq} \|\Delta^m_h f\|_p^q
\,dh/|h|^\ddd\biggr)^{1/q}
\]
or, in dimension $\ddd=1$,
%
%
\begin{equation}\label{ebeso}
|f|^s_{pq} = \biggl(2\int_0^1 h^{-1-sq} \|\Delta^m_h f\|_p^q
\,dh\biggr)^{1/q},
\end{equation}
where $\Delta^m_h$ denotes
the $m$th forward finite difference,
%
%
\begin{eqnarray}\label{ediff-sum}
\Delta^0_h f(x)&=&f(x),\nonumber\\
\Delta^m_h f(x)&=& [\Delta^{m-1}_h f(x+h)-\Delta
^{m-1}_hf(x)]\\
&=& \sum_{k=0}^m \pmatrix{m\cr k}(-1)^{m-k} f(x + k h).\nonumber
\end{eqnarray}
The Besov space $\Besov$ is the Banach space completion of $L_p(\bbR
^d)$ under
norm
%
%
\begin{equation} \label{ebesov-norm}
\|f\|^s_{pq} = \|f\|_p+|f|^s_{pq}.
\end{equation}
For $p = q = 2$, $\Besov$ coincides with the Sobolev space $\Sob$.

For fixed\vspace*{1pt} $\omega\in\Omega$, each of the kernel functions $g(\cdot,\omega)$ in
(\ref{ekernels}) is in $\Besov$ for all $p,q\ge1$ and some $s>0$,
and hence
each finite approximation of the form (\ref{eapprox-f}) lies in the same
$\Besov$. For example,\vspace*{1pt} the Gaussian kernel
of~(\ref{eGaussian})
(along with
its $d$-dimensional generalization)
satisfies $g_G(\cdot,\omega)\in\Besov$ for every\vadjust{\eject} $s<\infty$ and
$p, q \ge1$,
while in $\bbR^1$ the double-sided\vspace*{1pt} Laplace kernel of (\ref{eLaplace})
satisfies
$g_L(\cdot,\omega)\in\bbB^s_{pp}$ for $s<1 +1/p < 2$ for integer
$p$ and the
Haar wavelet of (\ref{eHaar}) is in $\Besov$ only for $s<1/p$. To simplify
proofs in Section \ref{ssgenBspq}, we will restrict attention to generating
functions $\GB$ on $\bbR^d$; these results may be extended to bounded domains
the Besov semi-norms defined in terms of differences on bounded domains in
Section 5.2.2 of \citet{Trie1992} may be used to extend these results.

We now provide conditions for LARK models to be in the same Besov space as
their generating functions.

\vspace*{3pt}\subsection{Convergence of LARK models in Besov spaces}\vspace*{3pt} \label{ssgenBspq}

\begin{theorem}\label{tbesov}
Fix $\GB\in\Besov(\bbR^d)$ for some $p, q \ge1$ and $s > 0$ and a
L{\'{e}}vy
measure $\nu$ on $\bbR\times\Omega$ with $\Omega=(\LamS\times
\bbR^d)$ of
translation-invariant product form $\nu(\dbdo) = \nubw(d\beta \, d\Scale) \,d\loc$ [here $\omega=(\Scale,\loc)$] for a $\sigma
$-finite measure
$\nubw(d\beta \, d\Scale)$ on $\bbR\times\LamS$ that satisfies the
integrability condition (\ref{el1bound}). Define a location-scale LARK
model $f(\cdot)$ on $\cX=\bbR^d$ by: $f(x) = \int_\Omega{\phi}(x,
\omega)
\Lmea( d \omega)$ where ${\phi}(x, \omega) \equiv\GB(\Scale
(x -
\loc))$ satisfies (\ref{enocomp}) for each fixed $x\in\cX$.
Then $f$
has the almost surely convergent series expression
%
%
\begin{equation}\label{ebesov-exp}
f(x) = \sum_j \GB\bigl(\Scale_j (x - \loc_j)\bigr) \beta_j
\end{equation}
and $f \in\Besov$ almost surely if $\nubw$ satisfies
%
%
\begin{subequation}\label{eLpBspq}
\begin{eqnarray}
\label{eLp}
\iint_\RS(1\wedge|\beta| |\Scale|^{- 1/p})
\nubw(d\beta \, d\Scale) &<&\infty,\\
\label{eBspq}
\iint_\RS(1\wedge|\beta| |\Scale|^{s - 1/p})
\nubw(d\beta \, d\Scale) &<& \infty.
\end{eqnarray}
\end{subequation}
\end{theorem}
\begin{pf}
Equation (\ref{enocomp}) ensures that the sum in (\ref{ebesov-exp}) will converge
almost~su\-rely for each fixed $x \in\cX$, with a finite number of terms
$|\GB(\Scale_j(x-\break\loc_j)) \beta_j | > 1$ and
infinitely
many, but absolutely summable, terms with\break $|\GB(\Scale_j(x-\loc_j)) \beta_j |\le1$.
The $L_p$ norm of $f$ satisfies the bound
\[
\|f\|_p \le\sum_j \bigl\|\GB\bigl(\Scale_j(\cdot-
\loc_j)\bigr) \bigr\|_p |\beta_j|
= \|\GB\|_p\sum_j |\Scale_j|^{-1/p} |\beta_j|
\]
by the triangle inequality and Proposition \ref{propdiff} in Appendix
\ref{aproofs}.
This is
finite almost surely by (\ref{eLp}) since $\GB\in\Besov\subset
L_p$. The
Besov semi-norm of $f$ is bounded by
\begin{eqnarray*}
|f|^s_{pq} &\leq& \sum_j |\beta_j|
\bigl|\GB\bigl(\Scale_j( x - \loc_j)\bigr)\bigr|^s_{pq} \\
& = &\sum_j |\beta_j| \biggl(\int_{|h|\le1} |h|^{-d -s q} \bigl\|\Delta^m_h
\GB\bigl(\Scale_j(\cdot- \loc_j)\bigr) \bigr\|^q_p \,dh \biggr)^{1/q} \\
& = &\sum_j |\beta_j| \biggl(\int_{|h|\le1|} |h|^{-d -s q} |\Scale
_j|^{-q/p} \|\Delta^m_{ \Scale_j h} \GB\|^q_p \,dh \biggr)^{1/q}
\end{eqnarray*}
by Proposition \ref{propdiff}; changing variables $h\mapsto t
= \Scale h$,
this is
%
%
\begin{equation}\label{ebsn}
= \sum_j |\beta_j| |\Scale_j|^{s - 1/p} \biggl(\int_{|\Scale
_j^{-1}t|\le1} |t|^{-d- s q} \|\Delta^m_{t} \GB\|^q_p \,dt \biggr)^{1/q}.
\end{equation}
The integral in (\ref{ebsn}) is bounded by
\begin{eqnarray*}
\int_{\bbR^d} |t|^{-d -s q} \|\Delta^m_{t} \GB\|^q_p \,dt
&=& \int_{|t|\le1} |t|^{-d- s q} \|\Delta^m_{t} \GB\|^q_p \,dt\\
&&{}+ \int_{|t|>1} |t|^{-d- s q} \|\Delta^m_{t} \GB\|^q_p \,dt.
\end{eqnarray*}
The first term is just $(|\GB|^s_{pq})^q$, and
(\ref{ediff-sum}) implies $\|\Delta^m_{t}\GB\|_p\le2^m\|\GB\|_p$, so
\begin{eqnarray*}
&\le&(|\GB|^s_{pq})^q
+ \int_{|t|>1} |t|^{-d-sq}(2^m\|\GB\|_p)^q \,dt\\
&=& (|\GB|^s_{pq})^q
+ \frac{\pi^{d/2} 2^{1+mq}} {\Gamma(d/2) s q}\|\GB\|^q_p
\\
&\le&(c \|\GB\|^s_{pq})^q
\end{eqnarray*}
for some $c<\infty$, so
%
%
\begin{equation} \label{eBesov-condition2}
|f|^s_{pq} \le c\|\GB\|^s_{pq} \sum_j |\beta_j| |\Scale_j|^{s - 1/p},
\end{equation}
which is almost surely finite by (\ref{eBspq}).
\end{pf}

Each of the kernels $g(\cdot,\omega)$ considered in the examples in
Sections \ref{ssim} and~\ref{sapps} may be shown to be in some
Besov space $\Besov$, and
each is\vspace*{1pt}
bounded by $\|\GB\|_\infty\le1$. Corollary \ref{cbes-egs}
establishes that
each of
our LARK models with a~L{\'{e}}vy measure that satisfies (\ref{el1bound})
is in
the same space $\Besov$ as its generating function.

\begin{cor} \label{cbes-egs}
Let $f(x)=\int{\phi}(x,\omega) \Lmea(d\omega)$ be a
one-dimensional LARK
model on a compact set $\cX\subset\bbR^1$, with product L{\'{e}}vy measure
$\nu(\dbdo)=\nu_\beta(d\beta) \pi_\lambda(d\lambda) \,d\loc$ on
$\bbR
\times\bbR^+ \times\cX$ satisfying (\ref{el1bound}) with Gamma probability
measure $\pi_\lambda(d\lambda) = \Ga(a_\lambda,b_\lambda)$ and
location-scale
generator ${\phi}(x,\omega)=\GB\uu$ with bounded $\GB\in\Besov$.
Then\vspace*{1pt} $f \in
\Besov$ almost surely if $\alpha_\lambda> 1/p$ for $p, q \ge1$ and
$s >
0$. In particular, if $a_\lambda\ge1$ then $f \in\Besov$ if
$\GB\in\Besov$ for all $p, q \ge1$ and $s > 0$.
\end{cor}

\begin{pf}
Equation (\ref{enocomp}) holds for bounded $\GB\in\Besov$ with L{\'{e}}vy
measures of
the form indicated; the conditions on $\alpha_{\lambda}$ ensure that also
$\int_{\bbR_+}\lambda^{-1/p}\pi_\lambda(d\lambda)<\infty$ and
$\int_{\bbR_+}\lambda^{s-1/p}\pi_\lambda(d\lambda)<\infty$, so
the bounds of
(\ref{eLpBspq}) hold.
\end{pf}

\subsection{Comparisons with Abramovich, Sapatinas and Silverman}
\label{ssass2k}
The sto\-chastic wavelet expansion of \citet
{AbraSapaSilv2000} may be
viewed as a LARK model using wavelet generator (\ref{eWave}), with
coefficients that, when conditioned on the scale parameters $\{a_j\}$,
have\vspace*{-1pt} independent Gaussian distributions $\{\beta_j\} \ind
\No(0, c a_j^{-\delta})$ with $\omega= (a, b) \in[a_0, \infty
)\times
[0,1)$ and $\nu_\omega(d \omega) \propto a ^{-\xi} \mathbf{1}_{\{a
\geq a_0\}} \,db \,da$
for some $c, \delta, \xi\geq0$, $\delta+ \xi> 0$ and $a_0 \geq1$.
The parameters $\delta$ and $\xi$ control the size and frequency of wavelet
coefficients and determine whether the expansion will have a well-defined
limit. For a~finite L\'{e}vy measure $\nu_\omega(d \omega)$ ($\xi>
1$), the
expansion will be in the corresponding Besov space of the generating wavelet
with probability one. For $\xi\leq1$, the Poisson mean is no longer
finite; however, \citet{AbraSapaSilv2000} provide conditions
on $\delta$
and $\xi$ so that $f$ falls in the corresponding Besov space of the
generating wavelet.

For ``simplicity of exposition,'' Abramovich, Sapatinas and Silverman
work with
functions of unit period [i.e., satisfying $\GB(x) = \GB(x+1)$] and regard
them as functions on the unit torus $\bbT$, the interval $[0,1]$ with the
endpoints identified. We now illustrate how the LARK theory may be used to
prove that the resulting expansion lies in $\Besov(\bbT)$
if the\vspace*{1pt}
generating function does.
%
%
The Besov sequence norms used by Abramovich, Sapatinas and Silverman
and others
are natural for the Gaussian distributions and discrete wavelet expansions
they study; we have found the (equivalent) function norms to be more
convenient for continuous wavelet expansions using non-Gaussian
($\alpha$-Stable, e.g.) distributions used for the coefficients
in our
expansions.
We follow \citeauthor{Niko1975} [(\citeyear{Niko1975}), Sections
1.1.1 and 4.3.5] in defining Besov norms on the torus by replacing the
$L_p$ norm on $\bbR$ with that over $\bbT$ in the definition of the Besov
semi-norm and norm [see (\ref{ebeso}), (\ref{ebesov-norm})], and in
denoting the
corresponding spaces by $L^*_p(\bbT)$ and $\Besper(\bbT)$,
respectively.

To simplify the proof, we will use the following lemma.
\begin{lem}\label{lASS}
Let $\pi_\Z(d\Z)$ denote the standard normal distribution on $\bbR
$, let
$\GB\in L_p^*(\bbT)$ with $p \ge1$
and let $r\in\{0,1\}$.
Then
\[
\iiint_\RRT
\bigl(1 \wedge|\Z\GB(u)^r | \lambda^{-a}\bigr) \lambda^{-b}
\pi_\Z(d\Z) \,d\lambda \,du < \infty
\]
for any $a\in\bbR$ if $b>1$, and for all $a>1-b$ if $b \le1$.
\end{lem}

The proof is given in Appendix \ref{alem1}.
\begin{theorem}\label{tbesov-nocomp}
Let $\GB\in\Besper(\bbT)$ for some $p, q \ge1$ and $s > 0$.
Let $\Lmea(d\omega)$ be a random field on $\Omega=[1,\infty)\times
\bbT$
with L{\'{e}}vy measure
%
%
\begin{equation}\label{eASS-nu}
\nu( d\beta \,d\lambda \,d\loc) =
\frac1{\sqrt{2\pi}}
\lambda^{{\delta}/{2}-\zeta}
e^{-\beta^2 \lambda^{\delta}/2}
\,d\beta \,d\lambda \,d\loc
\end{equation}
on $\bbR\,{\times}\,\Omega$ with $\delta, \zeta\,{\ge}\,0$. Then the LARK
model $f(x)\,{=}\,\int_\Omega\lambda^{1/2} \GB(\lambda(x\,{-}\,\loc))\Lmea(d\omega) $ has an
absolutely convergent
expansion 
%
%
\begin{equation} \label{eASS}
f(x) = \sum_j \beta_j\lambda_j^{1/2}\GB\bigl(\lambda_j(x-\loc_j)\bigr),\qquad 0\le
x<1,
\end{equation}
provided that $\frac{\delta-1}2 > 1-\zeta$ for $0 \le\zeta\le1$,
or for
any $\delta\ge0$ if $\zeta> 1$. Also $f(\cdot)\in\Besper(\bbT)$
almost surely for $\frac{\delta-1}2 > s+1-\zeta$ if $0 \le\zeta\le
1$ or
for any $\delta\ge0$ if $\zeta> 1$.\vspace*{-3pt}
\end{theorem}
\begin{pf}
The absolute convergence of (\ref{eASS}) for each $x$ will follow from
Proposition \ref{texist} if we can verify the conditions of (\ref
{enocomp}), that is,
finiteness of the integral
%
%
\begin{equation}\label{e18}
\iiint_\RRT\bigl(1 \wedge\bigl|\beta\lambda^{1/2} \GB\bigl(\lambda(x -
\loc)\bigr)\bigr|\bigr) \nu
(d\beta \,d\lambda
\,d\loc).
\end{equation}
Applying the change of variables
$\beta\mapsto\Z= \lambda^{\delta/2} \beta$,
%
%
\begin{equation}\label{eASSexist}
= \iiint_\RRT\bigl(1 \wedge|\Z| \lambda^{({1-\delta})/{2}} \bigl|\GB
\bigl(\lambda(x - \loc)\bigr)\bigr|\bigr) \lambda^{-\zeta} \pi_\Z( d\Z) \,d\lambda
\,d\loc,
\end{equation}
where $\pi_\Z(d\Z)$ is the standard normal
distribution. Since the term in parentheses is bounded by one,
(\ref{eASSexist}) is finite for all $\delta$ and $\GB$ if $\zeta>
1$. For
$0 \le\zeta\le1$, apply another change of variables $\loc\mapsto u =
\lambda(x - \loc)$ and apply periodicity
\[
= \int_\bbR\int
_1^\infty
\int_{ \lambda(x - 1)}^{\lambda x} \bigl(1 \wedge|\Z\GB(u)| \lambda
^{({1-\delta})/{2}}\bigr) \,du
\lambda^{-1-\zeta} \,d\lambda\pi_\Z( d\Z)
\]
which, due to periodicity, satisfies the bound
\begin{eqnarray*}
&\le&\int_\bbR\int_1^\infty\int_{0}^1 \bigl(1 \wedge|\Z\GB
(u)|\lambda^{({1-\delta})/{2}} \bigr) \,du
\lceil\lambda\rceil\lambda^{-1-\zeta} \,d\lambda\pi_\Z( d\Z)
\\[-3pt]
&\le&2 
\iiint_{\RRT} \bigl(1 \wedge|\Z\GB(u)|\lambda^{({1-\delta})/{2}}\bigr)
\,du \lambda^{-\zeta
} \,d\lambda
\pi_\Z( d\Z),
\end{eqnarray*}
where $\lceil\lambda\rceil$ denotes the least integer $\ge\lambda
$. By
Lemma \ref{lASS} this is finite for $0 \le\zeta\le1$ if $\frac
{\delta
-1}{2} >
1 - \zeta$ with $\GB\in\Besper$, 
so (\ref{enocomp}) holds and Proposition \ref{texist} ensures convergence.

The $L^*_p$ norms of the $m$th forward differences of a periodic
function
$\GB(\cdot)\in\Besper(\bbT)$ and their\vspace*{1pt} scaled translates $\GB
(\lambda( \cdot- \loc))$ for $\loc\in\bbT$ and positive scale
$\lambda\in[1,\infty)$ are related by
%
%
\begin{equation} \label{eqdiffnorm-per}
\bigl\| \Delta^m_h \GB\bigl(\lambda( \cdot- \loc)\bigr)\bigr\|^*_p
\le2^{1/p} \| \Delta^m_{\lambda h} \GB\|^*_p\vadjust{\eject}
\end{equation}
since, by a change of variables $x\mapsto u=\lambda(x-\loc)$,
\begin{eqnarray*}
&&
\bigl\| \Delta^{m}_h \GB\bigl(\lambda( \cdot- \loc)\bigr)\bigr\|^*_p
\\
&&\qquad= \lambda^{-1/p} \Biggl\{\int_{-\lambda\loc}^{\lambda(1 - \loc)} \Biggl|
\sum_{k=0}^m \pmatrix{m\cr k} (-1)^{m - k} \GB(u+k\lambda h) \Biggr|^p \,du
\Biggr\}^{1/p},
\end{eqnarray*}
which, again from periodicity, satisfies
\begin{eqnarray*}
&\le&\biggl(\frac{ \lceil\lambda\rceil}{\lambda}\biggr)^{1/p} \Biggl\{\int_{0}^{1}
\Biggl| \sum_{k=0}^m \pmatrix{m\cr k} (-1)^{m - k} \GB(u+k\lambda h) \Biggr|^p
\,du \Biggr\}^{1/p} \\
&=& \biggl(\frac{ \lceil\lambda\rceil}{\lambda}\biggr)^{1/p}
\| \Delta^m_{\lambda h} \GB\|^*_p,
\end{eqnarray*}
while $\lceil\lambda\rceil/\lambda\le2$.

The Besov semi-norm of $f$ is bounded by
%
%
\begin{eqnarray}\label{ebsn-per}
|f|^{s *}_{pq} & \leq& \sum_j |\beta_j| \lambda_j^{1/2} \bigl|\GB
{\bigl(\lambda_j(\cdot-\loc_j)\bigr)}\bigr|^{s *}_{pq}\nonumber\\
& =& \sum_j |\beta_j| \lambda_j^{1/2}
\biggl(\int_{|h|\le1} |h|^{-1 -s q} \bigl\|\Delta^m_h \GB\bigl(\lambda_j(\cdot-
\loc_j)\bigr) \bigr\|^{* q}_p \,dh \biggr)^{1/q} \nonumber\\
&\le&\sum_j |\beta_j| \lambda_j^{1/2}
\biggl(\frac{\lceil\lambda_j \rceil}{\lambda_j}\biggr)^{1/p}
\biggl(\int_{|h|\le1|} |h|^{-1 -s q} \|\Delta^m_{ \lambda_j h} \GB\|
^{* q}_p \,dh \biggr)^{1/q}\nonumber\\
& = & \sum_j |\beta_j| \lambda_j^{s + 1/2} \biggl(\frac{\lceil\lambda_j
\rceil}{\lambda_j}\biggr)^{1/p}\biggl(\int_{|t|\le\lambda_j} |t|^{-1- s q} \|
\Delta^m_{t} \GB\|^{* q}_p \,dt \biggr)^{1/q}.
\end{eqnarray}
The integral in (\ref{ebsn-per}) is bounded by
\begin{eqnarray*}
\int_{\bbR} |t|^{-1 -s q} \|\Delta^m_{t} \GB\|^{* q}_p \,dt
&=& \int_{|t|\le1} |t|^{-1- s q} \|\Delta^m_{t} \GB\|^{* q}_p \,dt\\
&&{} + \int_{|t|>1} |t|^{-1- s q} \|\Delta^m_{t} \GB\|^{* q}_p \,dt.
\end{eqnarray*}
The first term is just $(|\GB|^{s *}_{pq})^q$, and
(\ref{ediff-sum}) implies $\|\Delta^m_{t}\GB\|^*_p\le2^m\|\GB\|
^*_p$, so
\begin{eqnarray*}
&\le&(|\GB|^{s *}_{pq})^q
+ \int_{|t|>1} |t|^{-1-sq}(2^m\|\GB\|^*_p)^q \,dt\\
&=& (|\GB|^s_{pq})^q
+ \frac{2^{1+mq}} {s q}\|\GB\|^{* q}_p\\
&\le&(c \|\GB\|^{s *}_{pq})^q
\end{eqnarray*}
for some $c<\infty$, so
%
%
\begin{equation} \label{eBesov-condition-ASS}
|f|^{s *}_{pq} \le2c\|\GB\|^{s *}_{pq} \sum_j |\beta_j| \lambda
_j^{s + 1/2}
\end{equation}
is almost surely finite if and only if
\[
\iiint_{\RRT} (1 \wedge|\beta| \lambda^{s + 1/2})
\nu(d\beta \,d\lambda \,d\loc)
\]
is finite. Applying the change of variables $\beta\mapsto\Z=
\lambda^{\delta/2} \beta$,
\[
= \iiint_{\RRT}
\bigl(1 \wedge|\Z| \lambda^{s+ ({ 1-\delta})/{2}} \bigr) \lambda^{-\zeta
} \pi_\Z( d\Z) \,d\lambda \,d\loc
\]
is finite by Lemma \ref{lASS} for all $\delta\ge0$ if $\zeta> 1$
and for
$\frac{\delta- 1}{2} > s + 1 - \zeta$ if $0 \le\zeta\le1$. A~similar
argument
shows that the $L^*_p$ norm of $f$ satisfies a bound of the form
\[
\|f\|^*_p \le c \|\GB\|^*_p\sum_j |\beta_j| \lambda_j^{1/2}
\]
for some $c<\infty$. This is finite almost surely if
\begin{eqnarray*}
&&\iiint_{\RRT} (1 \wedge|\beta| \lambda^{1/2}) \nu(d\beta
\,d\lambda \,d\loc)\\
&&\qquad= \iiint_{\RRT} \bigl(1 \wedge|\Z| \lambda^{({1 - \delta})/{2}}\bigr)
\lambda^{-\zeta}
\pi_Z(d\Z)
\,d\lambda \,d\loc
\end{eqnarray*}
is finite, which follows from Lemma \ref{lASS} for all $\delta\ge0$ if
$\zeta>
1$ and, if $ \zeta\le1$, for~$\delta$ satisfying $\frac{\delta-
1}{2} > 1
- \zeta$ since $\GB\in\Besper\subset L^*_p$. Combining conditions,
the~$\Besper$
norm of $f$ is finite if $ \delta/2 - 1/2 > s + 1 - \zeta$ for $0 \le
\zeta
\le1$ and for all $\delta\ge0$ if~$\zeta> 1$.
\end{pf}

For L{\'{e}}vy measures $\nu(d\beta \,d\lambda \,d\loc)$ supported on
$\bbR\times
\bbN\times\bbT$ (i.e., for which $\lambda$ is almost-surely
integral) the
function $f(x)$ of (\ref{eASS}) would inherit periodicity from the generator
$\GB\uuj$ but, for the absolutely-continuous measure of (\ref
{eASS-nu}), it is
the definition of $f(x)$ as a function on $\bbT$ [as in \citet
{AbraSapaSilv2000},
equation (2)] that induces periodicity. The restriction to
$\lambda\ge1$ may be relaxed to the more natural $\lambda>0$ in the LARK
framework, but may require the use of compensation.

\subsection{Compensation}\label{ssBesovComp}
$\!\!\!\!$For L{\'{e}}vy measures satisfying only the local-$L_2$ bound of (\ref
{el2bound})
and not the local-$L_1$ bound of (\ref{el1bound}), we must use the definition
of $f(x)$ in (\ref{epois-comp}) and use (\ref{eexist}) to establish
conditions
that ensure $f$ will be\vspace*{1pt} well defined for $\GB\in\Besov$. We verify these
conditions for the existence of LARK models under symmetric $\alpha$-Stable
random fields.

\begin{theorem}\label{tstable}
For a Symmetric $\alpha$-Stable random field with L{\'{e}}vy measure
of the
form $\nu(\dbdo) = c_\alpha\alpha|\beta|^{-1 - \alpha} \,d\beta
\pi(d\Scale) \,d\loc$ on $\RS\times\Rd$ for \mbox{$0\,{<}\,\alpha\,{<}\,2$}, with
$\pi(d\Scale)$ a probability\vspace*{1pt} measure on $\LamS$ and $\GB\in
\Besov(\bbR^d)\cap L_1(\bbR^d)$ for $p, q \ge1$ and $s > 0$, the
conditions of (\ref{eexist}) for $f(x)$ to be well defined by Theorem
\ref{texist}
are satisfied for $1 < \alpha\le p$, $\alpha<2$ if $\E[|\Lambda|^{-1}]
< \infty$. For $\alpha= 1$, there is the additional requirement that
%
%
\begin{equation}
\label{eBesov-Cauchy}
\int_{\Rd} |{\GB(u) \log}|\GB(u)| | \,du < \infty.\vspace*{-6pt}
\end{equation}
%
\end{theorem}
\begin{pf}
Fix $x\in\cX$. By the affine change of variables of
\mbox{$\loc\mapsto u\equiv\Scale(x-\loc)$},
\begin{eqnarray*}
&&\iint_{[-1,1]^c \times\Omega} \bigl(1 \wedge|\beta{\phi
}(x,\omega)| \bigr)
\nu(\dbdo) \\
&&\qquad=2 {c_\alpha}{\alpha} \int_\LamS|\Lambda|^{-1} \pi(d\Lambda)
\iint_{[1, \infty) \times\Rd} \bigl(1 \wedge\beta|\GB(u)|\bigr)
\beta^{-1 - \alpha} \,d \beta \,du \\
&&\qquad = 2 {c_\alpha}{\alpha} \E|\Lambda|^{-1} \int_\Rd
\biggl\{\int_1^{|\GB(u)|^{-1}} \beta^{-\alpha} |\GB(u)| \,d\beta+ \int
_{|\GB(u)|^{-1}}^{\infty} \beta^{-1 -\alpha} \,d\beta\biggr\} \,du.
\end{eqnarray*}
For $1<\alpha<2$,
\[
= 2 {c_\alpha}{\alpha} \E|\Lambda|^{-1} \biggl\{\int_{\Rd} \frac
{|\GB(u)| - |\GB(u)|^\alpha}{\alpha-1} \,du + \int_{\Rd} \frac
{|\GB(u)|^\alpha} {\alpha}\,du \biggr\},
\]
which is finite for $1 < \alpha\le p$ since $\GB\in L_1$ and
$\GB\in\Besov\subset L_p$. For $\alpha= 1$,
\[
= 2{c_1}\E|\Lambda|^{-1} \biggl\{\int_{\Rd} -{|\GB(u)| \log}|\GB(u)| \,du +
\int_{\Rd} |\GB(u)|\,du \biggr\}.
\]
%
%
The first integral exists and is finite by (\ref{eBesov-Cauchy})
while the second is finite since $\GB\in L_1$.
%
Similarly, 
the integral in (\ref{eexb}) is
\begin{eqnarray*}
&&\iint_{[-1,1]\times\Omega}  \bigl(|\beta{\phi}(x,\omega) |
\wedge
|\beta{\phi}(x,\omega)|^2\bigr) \nu(\dbdo)\\
&&\qquad=2 c_\alpha\alpha\E|\lambda|^{-1}\biggl\{ \iint_{\ui\times\Rd}
\bigl(|\beta\GB(u) | \wedge|\beta\GB(u)|^2 \bigr)\beta^{-1 - \alpha}
\,d\beta \,du\biggr\}.
\end{eqnarray*}
The integral in braces
\[
\iint_{{[0,1\wedge|\GB(u)|^{-1}]\times\Rd}}
\beta^{1 - \alpha}
\GB(u)^2 \,d\beta \,du +
\iint_{{[1\wedge|\GB(u)|^{-1}, 1]\times\Rd}}
\beta^{ - \alpha} |\GB(u)| \,d\beta \,du
\]
is finite for $1 < \alpha\le p$, $\alpha<2$:
\begin{eqnarray*}
&\le&\int_\Rd\frac{|\GB(u)|^\alpha}{2 - \alpha}
+ \frac{|\GB(u)|^\alpha-|\GB(u)|}{\alpha-1} \,du \\
&\le&\frac{\|\GB\|_p^p}
{(2-\alpha)(\alpha-1)}< \infty,
\end{eqnarray*}
while for $\alpha= 1$,
\[
\le\int_\Rd\{|\GB(u)| + |{\GB(u) \log}|\GB(u)||\} \,du
<\infty
\]
by (\ref{eBesov-Cauchy}). Finally, (\ref{eexc}) holds because
%
\begin{eqnarray*}
&&
\iint_{\RO} (1 \wedge\beta^2) |{\phi}(x, \omega)| \nu
(\dbdo)\\
&&\qquad= \E|\Lambda|^{-1} c_\alpha\alpha
\iint_{\bbR\times\Rd}
(1 \wedge\beta^2) |\beta|^{-1 - \alpha} |\GB(u)| \,d\beta \,du \\
&&\qquad= \E|\Lambda|^{-1} \|\GB\|_1 c_\alpha\alpha
\int_\bbR(1 \wedge\beta^2) |\beta|^{-1 - \alpha} \,d\beta<
\infty.
\end{eqnarray*}
\upqed\end{pf}

All of the generator functions in the examples in Section \ref{ssim}
satisfy the
conditions of the theorem for the Cauchy random field ($\alpha= 1$),
so the
LARK models are well defined as $\eps\to0$ and for finite $\eps> 0$, the
approximations are in the same Besov space as $\GB$. We are able to show
that this also holds for Sobolev $\Sob$ spaces (which are equivalent to
$\bbB^s_{22}$) even when compensation is required, but this remains an open
question for $\Besov$ with general $p$ and~$q$.

\subsection{\texorpdfstring{Convergence in $\Sob$}{Convergence in W^s_2}}
\label{sSob}

\begin{theorem}\label{tsobo}
$\!\!\!$Let $\{{\phi}(x,\omega)\}$ be a location-scale family of the form
${\phi}(x,\omega)\!\equiv\G(\Ru(x-\loc))$ for $\omega
=(\loc,\Scale)$
with $\loc\in\bbR^d$ and nonsingular $d{\times}d$ matrix $\Scale
\in
\LamS$ for some function $\G(\cdot)\in\Sob$ with $s\ge0$. Let
$\nu$ be a
L{\'{e}}vy measure satisfying the condition
%
%
\begin{equation}\label{ethm2cond}
\iint_\RO|\Scale|^{-1} [1+\rho(\Scale)^{2s}] (1\wedge\beta^2)\nu(\dbdo) <\infty,
\end{equation}
where $\rho(\Scale)$ denotes the spectral radius (largest eigenvalue) of
$\Scale$. Recall
\begin{eqnarray*}
f(x) &\equiv&\iint_\RO{\phi}(x,\omega) [\beta-\cmp] \Np
(\dbdo)\\[-8pt]
\hspace*{-84pt}\mbox{(\ref{ef-def})}\hspace*{84pt}\\[-8pt]
&&{}+\iint_\RO
{\phi}(x,\omega) \cmp\Nt(\dbdo)
\end{eqnarray*}
and, for $\eps>0$, define
%
%
\begin{eqnarray}\label{efe-sum}
&&f_\eps(x)
\equiv\iint_{[-\eps,\eps]^c\times\Omega}
{\phi}(x,\omega) [\beta-\cmp] \Np(\dbdo)\nonumber\\
&&\qquad\quad{}+\iint_{[-\eps,\eps]^c\times\Omega}
{\phi}(x,\omega) \cmp\Nt(\dbdo)\\
&&\qquad=\mathop{\sum_{0\le j< J_{\eps}}}_{\eps< |\beta_j|}
{\phi}(x,\omega_j)\beta_j
-\iint_{[-\eps,\eps]^c\times\Omega}
{\phi}(x,\omega) \cmp\eta(\dbdo).\nonumber
\end{eqnarray}
Then $f_\eps(\cdot)\to f(\cdot)$ in $\Sob$ almost surely as $\eps
\to0$.
\end{theorem}
\begin{pf}
First, consider the case of compensator functions satisfying
$\cmp=\beta$ for all $|\beta|\le1$. Apply an affine change of
variables to see that ${\phi}(x,\omega)$ has Fourier transform (in $x$)
\[
\hat
{\phi}(\xi,\omega) = e^{i\xi\cdot\loc}|\Scale|^{-1}\hat\G
(\Rp^{-1}\xi).
\]
For
$0<\eps_1<\eps_2<1$ and $x\in\bbR^d$, set $\Delta(x) \equiv
f_{\eps_1}(x)-f_{\eps_2}(x)$ and let
$A\equiv\{\eps_1<|\beta|\le\eps_2\}\times\Omega$. Then
\[
\Delta(x)
=\mathop{\sum_{0\le j< J_{\eps_1}}}_{\eps_1 < |\beta_j|\le\eps_2}
{\phi}(x,\omega_j)\beta_j -
\iint_A 
{\phi}(x,\omega)\beta\nu(\dbdo)
\]
is a zero-mean random function of $x$ with Fourier
transform
\begin{eqnarray*}
\widehat\Delta(\xi)
&=&\mathop{\sum_{0\le j< J_{\eps_1}}}_{\eps_1 < |\beta_j|\le\eps_2}
e^{i\xi\cdot\loc_j}|\Scale_j|^{-1} \hat\G(\Rp_j^{-1}\xi)\beta
_j \\
&&{}-
\iint_A 
e^{i\xi\cdot\loc}|\Scale|^{-1} \hat\G(\Rp^{-1}\xi)\beta\nu
(\dbdo),
\end{eqnarray*}
a zero-mean $L_2$ random function of $\xi$ with second moment
%
%
\begin{equation}
\E|\widehat\Delta(\xi)|^2 = \iint_A
|\Scale|^{-2} |\hat\G(\Rp^{-1}\xi)|^2\beta^2
\nu(\dbdo). \label{el2Delt}
\end{equation}
Thus $\Delta(\cdot)$ has expected squared Sobolev norm
$\E\|{f_{\eps_1}-f_{\eps_2}}\|_\Sob^2$:
%
%
\begin{eqnarray}\label{enu-bound}
&=& (2\pi)^{-d}\iiint_{\bbR^d\times A} (1+|\xi|^2
)^s
|\Scale|^{-2} |\hat\G(\Rp^{-1}\xi)|^2\beta^2
\nu(\dbdo) \,d\xi\nonumber\\
&=& (2\pi)^{-d}\iiint_{\bbR^d\times A} (1+|\Rp\eta
|^2)^s
|\Scale|^{-1} |\hat\G(\eta)|^2\beta^2 \nu(\dbdo)
\,d\eta\nonumber\\
&\le&(2\pi)^{-d}\iiint_{\bbR^d\times A} (1+|\eta
|^2)^s
\bigl[\bigl(1+\rho(\Scale)\bigr)^{2s}\bigr] |\Scale|^{-1} |\hat\G(\eta
)|^2\beta^2
\nu(\dbdo) \,d\eta\hspace*{-25pt}\nonumber\\
&=& \|{G}\|_\Sob^2 \iint_{\{\eps_1<|\beta|\le\eps_2\}\times
\Omega}
\bigl[\bigl(1+\rho(\Scale)\bigr)^{2s}\bigr]
|\Scale|^{-1} \beta^2 \nu(\dbdo)\\
&\to&0 \qquad\mbox{as $\eps_1,\eps_2\to0$ by (\ref
{ethm2cond}),}\nonumber
\end{eqnarray}
so $\{f_{\eps_k}\}$ is a Cauchy sequence in $\Sob$ for any $\eps_k\to0$ and
\mbox{$\|{f-f_{\eps_k}}\|_\Sob\to0$}. Since $f_\eps$ is a finite linear
combination of
scaled translates of $\G\in\Sob$,\vadjust{\eject} each~$f_\eps$ (and hence~$f$)
lies in
$\Sob$ almost surely and Theorem \ref{tsobo} is proved for
compensator functions
satisfying $\cmp=\beta$ for $|\beta|<1$.

For an arbitrary bounded compensator $\cmp$ satisfying $|\beta-\cmp|
\le
c\beta^2$ for some $c>0$, (\ref{el2Delt}) has the additional
nonrandom term
\[
\biggl| \iint_A {e^{i\xi\cdot\loc}}
\frac{\hat\G(\Rp^{-1}\xi)} {| \Scale| }
\bigl(\beta-\cmp\bigr) \nu(\dbdo)\biggr|^2
\le c
\biggl( \iint_A
\frac{|\hat\G(\Rp^{-1}\xi)|} {| \Scale| }
\beta^2 \nu(\dbdo)\biggr)^2
\]
leading at most to an additional constant factor of
$[1\!+\!c\!\iint_\RO(1\!\wedge\!\beta^2) \nu(\dbdo)]$ in (\ref{enu-bound}),
leading as before to $\|{f-f_{\eps_k}}\|_\Sob\to0$ and completing
the proof.
\end{pf}

\begin{cor}\label{cPois}
If $\{{\phi}(x,\omega)\}$ is a location-scale family of the form
considered in
Theorem \ref{tsobo} and if a L{\'{e}}vy measure $\nu$ is of product
form $\nu
(\dbdo)
=\nu_\beta(d\beta) \pi_\omega(d\omega)$ for some $\sigma
$-finite measure
$\nu_\beta(d\beta)$ on $\bbR$ and probability measure $\pi_\omega
(\cdot)$
on $\Omega$ that for some $s\ge0$ satisfy
%
%
\begin{subequation}\label{el2sob}
\begin{eqnarray}
\label{el2sob-b}
\int_{\bbR} (1\wedge\beta^2) \nu_\beta(d\beta)
&<&\infty,
\\
\label{el2sob-w} 
\int_{\Omega} |\Scale|^{-1}\bigl(\bigl(1+ \rho(\Scale)\bigr)^{2s}\bigr)
\pi_\omega(d\omega) &<&\infty,
\end{eqnarray}
\end{subequation}
then $\nu(\dbdo)$ also satisfies (\ref{ethm2cond}) and hence
$f_\eps(\cdot)\to f(\cdot)$ in $\Sob$ almost surely as $\eps\to0$.
\end{cor}

For example, in one dimension, (\ref{el2sob-w}) is satisfied for all
$s>0$ if
$\Scale= \lambda$ has the $\chi_\nu$ distribution with $\nu>1$
degrees of
freedom, that is, if $\lambda^2\sim\Ga(\alpha_\lambda, \beta
_\lambda)$ with
$\alpha_\lambda>\half$. More generally, for any $m>0$ (\ref
{el2sob-w}) is
satisfied for all $s>0$ if $\lambda^m\sim\Ga(\alpha_\lambda, \beta
_\lambda)$
with $\alpha_\lambda>1/m$ or, for $m<0$, for $\alpha_\lambda>(1-2s)/m$.

Recall that the quantity $\eps$ introduced in the proof of Theorem
\ref{tsobo} and
the statement of Corollary \ref{cPois} is \textit{not} a model
parameter and has no
bearing on the Sobolov spaces to which the limiting function $f(\cdot)$
belongs; it is only a~tool used in proofs and implementations, to which we
now turn.

\section{Inference for LARK models}\label{sspreg} 

The LARK model introduced in Section~\ref{sintr} 
may now be summarized as
%
%
%
\begin{eqnarray}\label{elark-hier-all}
\label{lark-hier-a}
\E[Y(x)\mid\Lmea, \theta] &=& f(x) \equiv\int_\Omega
{\phi}(x,\omega)\Lmea(d\omega),\\
\label{lark-hier-b}
\Lmea\mid\theta&\sim&\Lv(\nu),\nonumber\\
\theta&\sim&\pi_\theta(d\theta)\nonumber
\end{eqnarray}
with implicit dependence of the L{\'{e}}vy measure $\nu(\dbdo)$ and
conditional
distribution for $Y(x)$ on a hyperparameter vector $\theta$. In all of our
examples, we take $\nu$ to be a product measure $\nu(\dbdo) =\nu
_\beta(\beta)
\,d\beta|\Omega| \pi_\omega(d\omega)$ satisfying the conditions of
Corollary \ref{cfe}, with $\pi_\omega(\cdot)$ a probability
measure on
$\Omega$,
$|\Omega|$ a~measure of the volume of $\Omega$, and $\nu_\beta
(\cdot)>0$ a
nonnegative density function on $\bbR$ satisfying $\int_\bbR(1\wedge
\beta
^2)\nu_\beta(\beta)\,d\beta<\infty$ [so $\nu$ satisfies (\ref
{el2bound})], for
which either (a) $\nu$ also satisfies (\ref{el1bound}) or (b)
$\nu_\beta(\beta)$ is even and $\cmp$ is odd in $\beta$. Thus, we
have the
representation
%
%
\begin{subequation}\label{elark-hier}
\begin{eqnarray}
\label{ehier-th}
\theta&\sim&\pi_\theta(d\theta),\\
\label{ehier-J}
J\mid\theta&\sim&\Po(\nu_\eps^+),\qquad
\nu^+_\eps\equiv\nu_\eps(\bbR\times\Omega),\\
\label{ehier-bw}
\{(\beta_j, \omega_j)\}\jJ\mid J,\theta
& \iid&\pi_\beta(\beta_j) \,d\beta_j
\pi_\omega(d\omega_j),\nonumber\\[-8pt]\\[-8pt]
\pi_\beta(\beta)&\equiv&\mathbf{1}_{\{|\beta|>\eps\}}
\nu_{\beta}(\beta)|\Omega|/\nu_\eps^+,\nonumber\\
\label{ehier-f}
Y_i\mid f & \ind& p\sY(y\mid f(x_i)) \,dy,\nonumber\\[-8pt]\\[-8pt]
f(x_i) &\equiv&\sum_{0\le j<J} {\phi}(x_i,\omega_j) \beta_j\nonumber
\end{eqnarray}
\end{subequation}
for sampling model $p\sY(\cdot\mid\mu)$ parametrized by $\mu$.

\subsection{\texorpdfstring{Examples of L{\'{e}}vy random fields}{Examples of Levy random fields}}\label{ssegRFs}
Motivated by the applications in Section~\ref{sapps}, we now focus on
LARK models
built on approximations to Gamma, symmetric Gamma and Symmetric
$\alpha$-Stable (in particular, Cauchy) L{\'{e}}vy random fields, and quantify
the approximation errors to facilitate the selection of $\eps$ and other
prior hyperparameters.

\subsubsection{Gamma LARK models}\label{sssinf-gam}
The Gamma random field of Section \ref{sssGRF} has $\nu_\beta
(d\beta)=
\gamma
\beta^{-1} e^{-\beta\eta}\mathbf{1}_{\{\beta>0\}} \,d\beta$ for
some constants
$\gamma>0$ and $\eta>0$. The parameter $\eta$ in (\ref{enu-gamma}) controls
both the Poisson rate of mass points $\{(\beta_j,\omega_j)\}$ of magnitude
$|\beta|>\eps$ and the probability distribution of those
magnitu\-des~$\{\beta_j\}$. To facilitate elicitation we disentangle those two
roles by
truncating at $|\beta\eta|\ge\eps$ (rather than $|\beta|\ge\eps
$); of course
the limit as $\eps\to0$ is the same. The distributions of $J$ and
$\{\beta_j\}$ are now given by
%
\begin{eqnarray*}
J &\sim& \Po(\nu^+_\eps),\qquad \nu^+_\eps= \gamma|\Omega|\Eone
(\eps),\\
\beta_j&\iid&\pi_\beta(\beta_j) \,d\beta_j ,\qquad
\pi_\beta(\beta_j) = \frac{\beta_j{}^{-1}
e^{- \beta_j\eta}}{\Eone( \eps)}\mathbf{1}_{\{\beta_j\eta>\eps\}},
\end{eqnarray*}
where the exponential integral function [\citet{AbraSteg1964},
pa\-ge~228] is
denoted as $\Eone(z)\equiv\int_z^\infty t^{-1}e^{-t} \,dt$. With this
truncation,
the expected square $L_2$ norm of the loss due to truncation
for any $\phi\in L_2(\Omega,\allowbreak|\Omega|\pi_\omega(d\omega)
)$, such as
$\phi(\omega)={\phi}(x,\omega)$, is
%
%
\begin{subequation}\label{el2err}
\begin{eqnarray}\label{el2errgam}
\E|\Lmea[\phi]-\Lmea_\eps[\phi]|^2 &=& \iint_{\bbR
\times\Omega}
%
\phi(\omega)^2 |\beta|^2 \mathbf{1}_{\{|\beta\eta|\le\eps\}}
\nu(\dbdo
)\nonumber\\
&=& \|\phi\|^2_2 \int_0^{\eps/\eta} \beta^2 \nu_\beta(\beta)
\,d\beta
\\
&=& \gamma\eta^{-2}\|\phi\|^2_2 [1-(1+\eps) e^{-\eps}
],\nonumber
\end{eqnarray}
\end{subequation}
showing the rate at which $\Lmea_\eps[\phi]\to\Lmea[\phi]$ in
$L_2$ as $\eps\to0$.
This is used in Section~\ref{sselic} to guide the elicitation of
hyperparameters.

\subsubsection{Symmetric Gamma LARK models}\label{sssymgamLARK}
The symmetric Gamma random field of Section \ref{sssSGRF} has L{\'
{e}}vy measure
$\nu_\beta(d\beta)= \gamma|\beta|^{-1} e^{-|\beta| \eta} \,d\beta$
for some constants $\gamma>0$ and $\eta>0$. Once again truncation at
$|\beta\eta|>\eps$ leads to
%
\begin{eqnarray*}
J &\sim& \Po(\nu^+_\eps),\qquad \nu^+_\eps= 2\gamma|\Omega|\Eone
(\eps)\\
\beta_j&\iid&\pi_\beta(\beta_j) \,d\beta_j,\qquad
\pi_\beta(\beta_j) = \frac{|\beta_j|^{-1}
e^{- |\beta_j|\eta}}{2 \Eone( \eps)}\mathbf{1}_{\{|\beta_j\eta
|>\eps\}}
\end{eqnarray*}
%
and expected squared discrepancy (used for elicitation)
%
\renewcommand{\theequation}{53b}
\begin{equation} \label{el2errsymgam}
\E|\Lmea[\phi]-\Lmea_\eps[\phi]|^2 =
2\gamma\eta^{-2} \|\phi\|^2_2 [1-(1+\eps) e^{-\eps}].
\end{equation}

\subsubsection{\texorpdfstring{Symmetric $\alpha$-Stable LARK models}{Symmetric alpha-Stable LARK models}}
\label{sssinf-SaS}

The \textsf{S$\alpha$S} L{\'{e}}vy random field of Section \ref{sssSaS}
has $\nu_\beta(d\beta)= \frac{\gd\alpha}{\pi}\Gamma(\alpha
)\sin\pat
|\beta|^{-\alpha-1} \,d\beta$
for some\vspace*{1pt} constants $\gd>0$ and $0<\alpha<2$. To facilitate
elicitation and
posterior inference, we write $\gd=\gamma\eta^{-\alpha}$ and (again)
truncate at $|\beta_j\eta|>\eps$.
This leads to
%
\begin{eqnarray*}
J &\sim& \Po(\nu^+_\eps),\qquad
\nu^+_\eps= \gamma|\Omega|\tpi\Gamma(\alpha)\sin\pat\eps
^{-\alpha}\\
\beta_j&\iid&\pi_\beta(\beta_j) \,d\beta_j,\qquad
\pi_\beta(\beta_j) = \frac{\alpha\eps^{\alpha}}{2 \eta
^\alpha} |\beta_j|^{-\alpha-1}
\mathbf{1}_{\{|\beta_j\eta|>\eps\}}
\end{eqnarray*}
with symmetric Pareto distributions for the coefficients $\{\beta_j\}
$. For the
Cauchy (\mbox{$\alpha=1$}), these simplify to
$\nu^+_\eps=2\gamma|\Omega|/(\pi\eps)$, with
\[
\pi_\beta(\beta_j) = {\frac{\eps}{2 \eta}} |\beta_j|^{-2}
\mathbf{1}_{\{|\beta_j\eta|>\eps\}}.
\]
Although the total variation $|\Lmea|$ is almost surely infinite, and
even
$|\Lmea-\Lmea_\eps|$ will be infinite for $\alpha\ge1$, still for
$\phi
\in L_2(\Omega,|\Omega|\pi_\omega(d\omega))$ the
expected squared
discrepancy is finite:
%
\renewcommand{\theequation}{53c}
\begin{eqnarray}\label{el2errSaS}
\E|\Lmea[\phi]-\Lmea_\eps[\phi]|^2 
&=& \iint_{\bbR\times\Omega}
\phi(\omega)^2 |\beta|^2 \mathbf{1}_{\{|\beta\eta|\le\eps\}}
\nu(\dbdo
)\nonumber\\[-8pt]\\[-8pt]
&=& 2\gamma\eta^{-2} \|\phi\|^2_2
\biggl[\frac{\Gamma(\alpha+1)} 
{\pi(2-\alpha)}\sin\pat\eps^{2-\alpha}\biggr]\nonumber
\end{eqnarray}
or $2\gamma\eta^{-2}\|\phi\|^2_2[\eps/\pi]$ for the Cauchy case
$\alpha=1$.

\subsection{Prior elicitation of hyperparameters}\label{sselic}

We now turn to the selection of $\eps>0$, the vector $\theta\in
\Theta$ of
(\ref{elark-hier}), and the L\'{e}vy measure $\nu(\dbdo)$. In each
of our
examples $\theta\equiv(\gamma,\eta)$ for rate parameters $\gamma$ and
$\eta$ governing the frequency and magnitude of coefficients $\{\beta
_j\}$,
respectively, and the expected squared truncation error for
$\Lmea_\eps[\phi]$ for $\phi(\omega)={\phi}(x,\omega)$ is of
the form
$\E|\Lmea[{\phi}(x,\cdot)]-\Lmea_\eps[{\phi}(x,\cdot)]
|^2 =\gamma\eta^{-2}
\|{\phi}(x,\cdot)\|_2^2 c(\eps)$ for some $c(\eps)>0$ with
$c(\eps)\to0$ as
$\eps\to0$ [see (\ref{el2err})].

We choose prior distributions to attain three goals: (1) desired range of
number $J$ of terms in the stochastic expansion; (2) desired range of
coefficient magnitudes $\{\beta_j\}$; and (3) tolerable expected truncation
error. We first select a L\'{e}vy family (Gamma, $\alpha$-Stable,
etc.) to
meet the needs of a particular problem for symmetry or positivity,
sharp or
heavy tails, etc. Each of our L\'{e}vy measures is of the product form
$\nu(d\beta \,d\omega) =\nu_\beta(d\beta) \pi_\omega(d\omega)$
considered
in Theorem \ref{tsobo} and Corollary \ref{cPois}, with location,
scale, and perhaps other
location-specific (and hence adaptive) attributes encoded in
$\omega\in\Omega$ in problem-specific ways.

Hyperparameters in the L\'{e}vy measure $\nu_\beta(d\beta)$ govern sparseness
for LARK models, that is, the number $J$ of terms in the stochastic expansion.
In each LARK model, $J$ has a Poisson distribution with mean
proportional to
$\gamma$. The coefficient of variation under the Poisson distribution falls
to zero as the mean increases, overstating the prior certainty for large
values of $\E J$.
To ameliorate this, we introduce an additional layer of hierarchy by
placing a
Gamma prior distribution on the parameter $\gamma\sim\Ga(a_\gamma,
b_\gamma)$, leading to the overdispersed negative binomial prior distribution
for $J\sim\NB(a_J, p_J)$. The parameter $\eta$ governs the scale of the
coefficients $\{\beta_j\}$, and hence the range of the regression function
$f(\cdot)$. We employ a Gamma distribution for the scale parameter
$\eta^{-1}\sim\Ga(a_\eta,b_\eta)$. Together the hyperparameters
$\eps$, $a_\gamma$, $b_\gamma$, $a_\eta$, $b_\eta$ determine the prior distributions
for $J$, for the coefficients $\{\beta_j\}$ (and hence the range of
$f(\cdot)$), and for the expected mean-square truncation error.
We select values for these five parameters to meet five criteria:
attain two
specified quantiles (such as a central 99\% interval) for each of $J$ and
$\{\beta_j\}$, and a specified bound on the expected truncation error
$\E\gamma\eta^{-2} \|{\phi}(x,\cdot)\|_2^2 c(\eps)$. Typically
this involves an
iterative numerical solution.

As a default choice, we take $\pi(d\omega)= \pi_\loc(d\loc)
\pi_\lambda(d\lambda)$ to be the product of the uniform distribution for
locations $\loc\sim\Un(\cX)$ and a Gamma distribution for inverse (distance)
scale parameters $\lambda\sim\Ga(a_\lambda, b_\lambda)$. The shape
and rate
hyperparameters $a_\lambda$ and $b_\lambda$ govern the range of
probable values
for the location-specific inverse scale parameters $\{\lambda_j\}$ and hence
for the smoothness of $f(x)$, similar to how bandwidth selection governs
smoothness in other kernel methods. A kernel at $\omega_j=(\loc
_j,\lambda_j)$
will represent a feature located at $\loc_j$ of width $1/\lambda_j$,
so large
values of $\lambda_j$ are needed to fit a very ``spiky'' part of a curve,
while a smoother part of a curve may be fit most parsimoniously using small
values of $\lambda_j$. The prior distribution for $\lambda_j$ must\vadjust{\eject}
support an adequate range of values in order to fit a spatially inhomogeneous curve.
Values of $a_\lambda>1$ will ensure $\E[\lambda]<\infty$ and a finite
covariance function; we choose $(a_\lambda,b_\lambda)$ to attain two specified
quantiles, such as a central $99\%$ interval.

\subsection{Posterior inference}\label{sspost}

The joint posterior density of all parameters under the LARK model of
(\ref{elark-hier}), given observations $\bY=\{Y_i\}$, is
%
%
\renewcommand{\theequation}{\arabic{equation}}
\setcounter{equation}{53}
\begin{eqnarray} \label{epostall}
&&p(\gamma, \eta, J, \bbet, \bw\mid\bY)\nonumber\\
&&\qquad\propto\pi_\gamma(\gamma) \pi_\eta(\eta) 
\frac{\exp[-\nu_{\eps}(\bbR\times\Omega)]}{J!}
\\
&&\qquad\quad{} \times\biggl\{ \prod_{0\le j<J} \nu_\eps(\beta_j, \omega_j) \biggr\}
\biggl\{\prod\iI p_Y\biggl(Y_i \Bigm| 
\sum_{0\le j<J} {\phi}(x_i, \omega_j)
\beta_j\biggr)\biggr\}.\nonumber 
\end{eqnarray}
The posterior (and full conditional) distributions of the parameters
are not
available in closed form. Since some of our parameters ($\bbet$ and
$\bw$)
have varying dimension, some form of trans-dimensional Markov chain Monte
Carlo, such as a reversible jump (RJ-MCMC) algorithm [\citet{Gree1995},
\citet{WolpIcksHans2003}, \citet{Siss2005}] must be used to
provide samples from
(\ref{epostall}) for posterior inference. See Appendix \ref
{arjmcmc} for a sketch
of the
RJ-MCMC algorithm.

\section{Relation of LARK to other models}\label{sother}

\subsection{Gaussian processes or random fields} \label{ssGRF}
For any positive Borel measure $\Sigma(d\omega)$ on a complete separable
metric space $\Omega$, there exists a Gaussian random measure $\cZ
(d\omega)$
on $\Omega$ that assigns to disjoint Borel sets $A_i\subset\Omega$
of finite
measure $\Sigma(A_i)<\infty$ independent mean-zero Gaussian random variables
$\cZ(A_i)\sim\No(0, \Sigma(A_i))$ of variance
$\E\cZ(A_i)^2=\Sigma(A_i)$. For any kernel function $g$ on $\cX
\times\Omega$
with ${\phi}(x,\cdot)\in L_2(\Omega,\Sigma(d\omega) ) $
for each
$x\in\cX$, this induces a mean-zero Gaussian random field through the Wiener
stochastic integral
\[
f(x) = \int_\Omega{\phi}(x,\omega) \cZ(d\omega)
\]
with covariance $C(x,y) = \E[ f(x) f(y) ] = \int_\Omega{\phi}(x,\omega){\phi}(y,\omega) \Sigma(d\omega)$.
The Gaussian random measure $\cZ(d\omega)$ is the special case of a L\'{e}vy random measure~$\Lmea(d\omega)$ defined
earlier in (\ref{elevy-full}) with $\delta(d\omega)\equiv0$ and $\nu(\dbdo)\equiv0$.

A wide variety of Gaussian processes are available in this form. For
example, those with stationary covariance $C(x,y)=c(x-y)$ may be
written in
the above form if the spectral measure has a density function $\hat
c(\omega)
=\int_\cX e^{-i\omega\cdot x} c(x) \,dx$ whose square root is Lebesgue
integrable, for example, the Mat\'{e}rn class [\citet{Stei1999},
page 31] in $\bbR^d$
with smoothness parameter $\nu>d/2$. The Gaussian random field model above
may also be obtained as the limit as $\alpha\to2$ of the symmetric
$\alpha$-Stable LARK models considered herein, providing an alternative
method for inference that avoids the need for large matrix inversions. To
maintain a unified computational approach, we have limited our
attention in
this article to LARK models with pure-jump L\'{e}vy random measures,
that is,
$\Sigma(\cdot) \equiv0$.

\subsection{Compound Poissons and mixtures of Gaussian random fields}
Mixtures of Gaussian random fields may be constructed as LARK models with
L\'{e}vy measure of the form
%
%
\begin{equation}\label{enormal}
\nu(\dbdo) = (2\pi\sigma_\omega^2)^{-1/2}
e^{-\beta^2/2\sigma_\omega^2} \,d\beta\nu_\omega(d\omega)
\end{equation}
leading to mean functions of the form $f(x_i) = \sum_{0\le j<J}
{\phi}(x,\omega_j) \beta_j$ with nor\-mally-distributed coefficients
$\beta_j|\omega\sim\No(\mu_\omega, \sigma_\omega^2)$. For
finite measures
$\nu_\omega$, the expansion has a Poisson-distributed number of
terms, hence,
is a Poisson mixture of Gaussian processes (or for hierarchical models
with a
Gamma distributed Poisson mean, a~negative binomial mixture of Gaussian
processes). In Section \ref{ssass2k}, we showed that the stochastic wavelet
expansion of \citet{AbraSapaSilv2000}, an example of (\ref
{enormal}), may
be viewed as a LARK model. \citet{ChuClydLian2009} extend the compound
Poisson (or LARK with finite $\nu$) model to include mixtures of normals
distributions for $\beta_\omega$ and develop methods for Bayesian inference
for such OverComplete Wavelet expansions (OCW); we compare the OCW
method to
other LARK models in the simulation study of Section \ref{ssim}.

For automatic curve fitting using splines and wavelets, Denison et al.
[(\citeyear{DeniHolmMallSmit2002}), Chapter 3] used a similar
hierarchical model
with common \mbox{$\sigma_\omega\equiv\sigma$}, but truncated the
(Poisson-distributed) number of terms in the basis expansions at some fixed
upper bound $J_u$. Taking $J_{u} \to\infty$ leads to the Gaussian LARK
model of (\ref{enormal}) with a common variance. Gaussian processes have
sharp tails, of course, leading to concerns about robustness when they are
used as prior distributions in problems with likelihood functions that fall
off more slowly.
Specifying variances for Gaussian prior distributions is nontrivial, with
large ``noninformative'' choices leading to the so-called Lindley paradox.
Denison et al. recommend an inverse Gamma prior on~$\sigma^2$
to avoid this well-known problem. This leads to a multivariate
Student~$t$ distribution on the expansion coefficients and, since the prior
now has bounded influence, provides robustness. The limiting model (as
$J_u\to\infty$) may be viewed as a~mixture of L\'{e}vy random fields.

Rather than using a multivariate Student $t$ for the coefficients, one might
use ``ridge'' priors and model\vspace*{1pt} the uncertain function $f(\cdot) =\sum
_{0\le
j< J}\beta_j {\phi}(\cdot;\omega_j)$ as the sum of a Poisson
(or negative
binomial)-distributed number $J$ of kernel functions ${\phi}(\cdot;
\omega_j)$
with coefficients $\beta_j \iid\Ca(0,\tau)$ drawn from a centered Cauchy
distributions with scale $\tau$. To accommodate rough functions
$f(\cdot)$, one must be willing to consider large numbers of terms, most of which will
have small coefficients---under these priors, one must consider large
$\E J$
and small $\tau$. But how small? And what happens if $\tau$ is made a bit
smaller and $\E J$ a bit larger? As $\tau\to0$, if one scales the expected
number $\E J$ of terms (as a function of $\tau$) properly, this model
converges to a LARK model with infinite L\'{e}vy measure (and so is
\textit{not}
sensitive to the cut-off $\varepsilon$, which merely quantifies how
close is
this approximation). If $\E J$ is not scaled properly to converge to a LARK
model, the limiting results may depend critically on arbitrary and
unintentional choices.

This may\vspace*{1pt} be implemented explicitly in LARK form by placing independent
$\Ga(\alpha/2, \eps/2)$ prior distributions on $\sigma^{-2}_{\omega
}$ in
(\ref{enormal}) to achieve independent univariate Student $t_\alpha(0,
\eps)$
distributions for the coefficients $\{\beta_j\}$ and (approximately,
as the
parameter $\eps\to0$) the heavy-tailed Symmetric $\alpha$-Stable
process for
$f(x)$ of Sections \ref{sssSaS} and \ref{sssinf-SaS} [this also
illustrates that truncating
the support of $\beta_\omega$ is not the only way to construct suitable
approximating sequences of finite L{\'{e}}vy measures $\nu_\eps(\dbdo)
\Rightarrow\nu(\dbdo)$ for which the integrals in
(\ref{epoi-comp-res}) 
converge]. An important feature of our infinitely divisible
construction (in
contrast to a compound Poisson approach from other distributional families)
is that in each case, as $\eps\to0$ 
the approximating model converges to one with a well-defined prior (with
infinite L{\'{e}}vy measure) and a proper posterior distribution.

\subsection{Finite dimensional frames}
LARK may be viewed as a limit of Bayesian variable selection methods with
finite frames or dictionaries. \citet{WolfGodsNg2004} consider frames
based on discretizing $\Omega$ as a fine grid with $|G|$ elements. They
place i.i.d. prior distributions $\pi_G(\beta) \,d\beta$ on the nonzero
coefficients and i.i.d. Bernoulli kernel inclusion indicators with inclusion
probability $\rho_G$. If $|G|\rho_G \pi_G(\beta)\to\nu(\beta)$ as
$|G|\to\infty$, then the result converges to a LARK model on the
infinite-dimensional frame.
The representation in \citet{WolfGodsNg2004} uses a point mass
at zero to
provide sparsity. Similarly, one may view the prior distributions in LARK
under the $\varepsilon$-truncation approach as assigning zero mass to a
neighborhood around zero, also leading to sparse representations. One
benefit of LARK is its provision of a formal method for coherent prior
specification for continuous dictionaries; a second is its provision of a
proper prior specification in the limit as $\eps\to0$, ensuring insensitivity
to the choice of $\eps$.

Standard stochastic search algorithms using finite-dimensional frames may
exhibit poor mixing when the correlations between grid elements tend to~$\pm1$.
To illustrate, suppose that two possible kernel parameters
$\omega_0$ and $\omega_1$ are close in parameter space, leading to
two highly
correlated columns in the design matrix. In addition, assume that inclusion
of either column leads to nearly-maximal likelihood. With the standard
one-at-a-time deletion or addition moves in many stochastic search
algorithms, to move from a model\vadjust{\eject} including a kernel indexed by $\omega
_0$ to
one indexed by $\omega_1$ would require an extremely unlikely deletion
followed by an addition (or unlikely addition followed by a~deletion). LARK
avoids this difficulty by allowing the continuous parameter~$\omega$ indexing
dictionary elements to move incrementally from~$\omega_0$ to~$\omega_1$
by a~series of update steps, avoiding some of the poor mixing problems associated
with highly correlated frame elements in a fine-grid based
method.

\subsection{Dirichlet processes}\label{ssdir}
The Dirichlet process [Ferguson (\citeyear{Ferg1973,Ferg1974}),
Antoniak (\citeyear{Anto1974})] has received
widespread use as a prior distribution on probability distribution functions.
Its popularity is due in large part to its analytic tractability in many
problems; simulation is straightforward, and Bayesian MCMC inference methods
are available [\citet{Esco1994}, \citet{MacE1994},
\citet{EscoWest1995}, \citet{MacE1998},
\citet{MullQuin2004}]. \citet{LianMukhWest2006} consider nonlinear
regression and classification models $\E[Y_i\mid X_i]= f(X_i)$ for data
$\{(Y_i,X_i)\}$ using kernel expansions of the form
%
%
\begin{equation}\label{eqLMW}
f(x) = \int k(x, u) \gamma(du) = \int k(x, u) w(u) F(du)
\end{equation}
with random signed measure $\gamma(du)$ expressed as the integral of a weight
function $w(u)$ with respect to a probability distribution $F$, modeled
as a~Dirichlet process $F\sim\DP(F_0, \alpha)$ with base measure $F_0$
and scale $\alpha>0$. If observed points $\{X_i\}$ are viewed as a random sample from
$F$, then updating the posterior for $F$ solely on the basis of the observed
$\{X_i\}$ would lead in the limit as $\alpha\to0$ to a degenerate posterior
for $F$ concentrated at the empirical distribution for $X$, justifying the
finite-dimensional expansion
\[
f(x) = \sum_{i=1}^n k(x, x_i) w(x_i)
\]
with kernels evaluated only at the observed data locations. The generalized
$g$-prior of \citet{West2003} for the coefficients $\{
w_i=w(x_i)\}$ leads to
dependent Cauchy distributions for the $\{f(x_i)\}$.
This approach (like the SVM, RVM and related approaches) has as many
coefficients as there are data points, 
but avoids over-fitting through shrinkage. Asymptotic properties of $f(x)$
as $n\to\infty$ are difficult to study in the absence of a limiting structure
such as that provided by LARK.

The Dirichlet measure $F(du)$ does not assign independent random
variables to
disjoint sets and so (\ref{eqLMW}) is not a LARK model, but it can be
constructed from one. In fact it is exactly the \textit{normalized} LARK model
%
\begin{eqnarray}\label{enormgam}
f(x) &=& \int_{\Omega} k(x, u) w(u) \Lmea(d u) \big/ \Lmea(\Omega)
\nonumber\\[-8pt]\\[-8pt]
&=& \sum_j k(x, u_j) w_j \beta_j \big/ \beta_+\nonumber
\end{eqnarray}
with $F(du)=\Lmea(du)/\Lmea(\Omega)$ for a Gamma random field
$\Lmea(d u)$ with infinite L\'{e}vy measure
\[
\nu(d\beta \,d u) = \alpha\beta^{-1} e^{-\beta}
\mathbf{1}_{\{\beta>0\}}\,d\beta F_0(du),
\]
where $\beta_+:=\sum\beta_j$ [note that $w(u)$ could be absorbed into
$k(x,u)$].

Well-known disadvantages of Dirichlet process models include their
inflexibility (the single parameter $\alpha$ determines the prior dispersion
\textit{everywhere}, precluding prior specifications with more
uncertainty in
some regions than in others), their discreteness, and the limited variability
of the masses assigned to the countably-many support points. The normalized
Gamma representation (\ref{enormgam}) of DP's offers the opportunity to
overcome some of these disadvantages---for example, the Gamma process
may be
given a variable rate parameter $b(u)$ by taking
\[
\nu(d\beta \,du) = \beta^{-1} e^{-b(u) \beta} \mathbf{1}_{\{\beta
>0\}}\,d\beta F_0(du)
\]
leading to a precision that can vary with location $u\in\Omega$, or
the Gamma
random field may be replaced with another nonnegative L{\'{e}}vy
random field
with wider dispersion, such as the fully-skewed Stable process of index
$\alpha<1$. Other nonnegative L{\'{e}}vy random fields are beginning to
be used
in machine learning [\citet{Jord2010}] and other fields.

\section{Simulation study}\label{ssim}

We now turn our attention to simulated and real examples to illustrate the
performance of LARK models in practice. We conducted a simulation study
using four spatially varying functions introduced by \citet
{DonoJohn1994}
that are now standard in the wavelet literature: Blocks, Bumps, Doppler and
Heavysine. Data were generated for each test function by adding independent
Gaussian random noise $\No(0,\sigma^2)$ to the true target function
$f(\cdot)$ at $n=1024$ equally-spaced points on $\cX=[0,10]$. As in
\citet{AbraSapaSilv1998}, the value of $\sigma$ was chosen to
attain a root
signal-to-noise ratio (RSNR) of $\sqrt{\int_\cX(f(x)-\bar{f})^2 \,dx
/\sigma^2} = 7.0$, where $\bar{f}\equiv\frac1{|\cX|}\int_\cX
f(x) \,dx$.
Each target function $f(\cdot)$ has a range of approximately $0\le
f(x)\le
25$. For each function, we generated $100$ replicate data sets to evaluate
the performance of LARK and other methods on the basis of mean squared error
%
\begin{equation}\label{emse}
\mathrm{MSE} \equiv n^{-1}\sum_{i=1}^n \bigl(\widehat
{f(x_i)}-f(x_i)\bigr)^2.
\end{equation}

\subsection{Hyperparameters}\label{sshyper}

In Table \ref{ttestf}, we report the kernel functions used for the
four simulation
%
\begin{table}
\tablewidth=250pt
\caption{Kernel functions used for four test functions}
\label{ttestf}
\begin{tabular*}{\tablewidth}{@{\extracolsep{\fill}}lc@{}}
\hline
\textbf{Test function} & \textbf{Kernel} $\bolds{{\phi}(x_i; \loc_j, \lambda_j)}$ \\
\hline
Blocks& $ \mathbf{1}_{\{0 < \lambda_j (x_i - \loc_j) \le1\}} $ \\[2pt]
Bumps& $ e^{-\lambda_j |x_i - \loc_j|} $ \\[2pt]
Doppler& $ e^{-0.5\lambda_j^2 (x_i - \loc_j)^2} $\\[2pt]
Heavysine& $ e^{-0.5\lambda_j^2 (x_i - \loc_j)^2}\mathbf{1}_{\{|x_i
- \loc_j| < 2.0\}} $\\
\hline
\end{tabular*}
\end{table}
examples, chosen to illustrate the flexibility of LARK to use a wide
range of
kernels that may be adapted to anticipated features (smoothness, spikiness,
jumps, curvature, covariation, etc.) of applications. In each case,
we take
$\Omega=[0,10]\times\bbR_+$ (and $|\Omega|=10$), with elements denoted
$\omega=(\loc,\lambda)$, comprising a location parameter $\loc\in
\cX=[0,10]$
and a~shape parameter $\lambda>0$. As described in Section \ref
{sselic}, we take
$\{\loc_j\}\iid\Un(\Omega)$ and $\{\lambda_j\}\iid\Ga(a_\lambda
,b_\lambda)$ with
$a_\lambda$, $b_\lambda$ chosen (see Table \ref{thyp}) to achieve a
95\% prior
interval of $[0.20,20.0]$ for $\lambda$ to attain dilated kernels
covering from
half a~percent up to fifty percent of $\cX$.


%
\begin{table}
\tablewidth=320pt
\caption{Hyperparameters used in examples of Section
\protect\ref{sshyper}}\label{thyp}
\begin{tabular*}{\tablewidth}{@{\extracolsep{\fill}}lccccccc@{}}
\hline
\textbf{L\'{e}vy measure} &$\bolds{\eps}$ &$\bolds{a_\gamma}$&$\bolds{b_\gamma}$
&$\bolds{a_\eta}$ &$\bolds{b_\eta}$
&$\bolds{a_\lambda}$&$\bolds{b_\lambda}$\\
\hline
Symmetric Gamma &0.0041&2.53 &\hphantom{0}6.45
&13.01 &0.71
&1.117 &0.1965\\
[3pt]
Cauchy &0.0029&2.53 &14.2\hphantom{0}
&\hphantom{0}0.50 &1.00
&1.117 &0.1965\\
\hline
\end{tabular*}
\end{table}

Our choice of the remaining hyperparameters was guided by three objectives:
to achieve a 95\% prior predictive interval of $[5,100]$ for $J$, to achieve
a 95\% prior predictive interval of $[-25,25]$ for the $\{\beta_j\}$,
and to
achieve a~limit on the mean squared truncation error of $
\|\Lmea[\phi]-\Lmea_\eps[\phi]\|_2 = (\E|\Lmea[\phi
]-\Lmea_\eps
[\phi]|^2 )^{1/2} \le0.05 \cdot\|\phi\|_2$ (see Section
\ref{sselic}).
While these objectives could be met for the LARK model with symmetric Gamma
prior with the values given in Table \ref{thyp}, they are not quite
attainable for
the Cauchy model---the competing goals of an extremely wide
distribution for
the $\{\beta_j\}$ and a low mean squared truncation error cannot be
reconciled. Upon relaxing the prior predictive distribution requirement on
$\{\beta_j\}$ to a
99.9\% interval of $[-33, 33]$, adequate for this problem with a flat
Pareto-tailed distribution for $\{\beta_i\}$, the remaining objectives for
the distribution of $J$ and the mean square truncation error were attained
using the values given in Table \ref{thyp}. See Figure~\ref{fRealize}, Appendix \ref{aegs}
for realizations
from the prior distribution.

%

\subsubsection{Performance}\label{sssperf}
$\!\!\!\!\!$We compared LARK with two of the best wavelet~me\-thods currently
available for
inhomogeneous function estimation using~over\-complete representations: the
empirical Bayes approach (``EBayesThresh'') of Johnstone and
Silverman (\citeyear{JohnSilv2004b,JohnSilv2005b,JohnSilv2005a}) using
translational-invariant wavelets,
and the continuous over-complete wavelet (``OCW'') approach of
\citet{ChuClydLian2009} based on the stochastic wavelet
expansions of
\citet{AbraSapaSilv2000}. We replicated the results of
\citet{JohnSilv2005a} under the beta-Laplace prior using their
\texttt{R} package
\texttt{EBayesThresh}
[\citet{JohnSilv2005b}] with \citeauthor{Daub1988}'
``least asymmetric'' (\texttt{la8}) wavelets [see Section 4
of \citet{Daub1988} or Section 6.4 of \citet{Daub1992}].
OCW uses the same
\texttt{la8} wavelet as EBayesThresh except for the Blocks example, where
both LARK and OCW use the Haar wavelet. The OCW method may be viewed as a
special case of LARK with a finite nonseparable L{\'{e}}vy measure, where
coefficients~$\beta_j$ have independent Laplace distributions
conditional on
scale parameters~$\lambda_j$, which in turn have truncated Pareto
distributions. As in LARK, OCW assigns independent uniform locations,
with a negative binomial distribution for the number of terms in the expansion.

%
%
\begin{table} 
\caption{Average and (standard errors) over $100$ replications of mean square
errors of the four test functions using the L{\'{e}}vy Adaptive Regression
Kernels (LARK) using the symmetric Gamma
and Cauchy priors, the OCW approach using a Laplace prior
[Chu, Clyde and Liang (\protect\citeyear{ChuClydLian2009})],
and the EBayesThresh approach using
a Laplace prior [Johnstone and Silverman (\protect\citeyear{JohnSilv2005b})]}
\label{tmse}
\begin{tabular*}{\tablewidth}{@{\extracolsep{\fill}}lcclc@{}}
\hline
\textbf{Method} & \textbf{Blocks} & \textbf{Bumps} & \multicolumn{1}{c}{\textbf{HeavySine}}
& \textbf{Doppler} \\
\hline
LARK-Gamma& 0.030 (0.0013) & 0.111 (0.0019) & 0.038 (0.0010) & 0.152
(0.0030)\\
LARK-Cauchy& 0.026 (0.0011) & 0.105 (0.0017) & 0.036 (0.0010) & 0.157
(0.0028)\\
OCW & 0.060 (0.0023) & 0.285 (0.0025) & 0.082 (0.0010) & 0.152
(0.0019)\\
EBayesThresh & 0.096 (0.0013) & 0.307 (0.0032) & 0.118 (0.00098)
& 0.202 (0.0027) \\
\hline
\end{tabular*}
\end{table}

%
\begin{figure}[t!]

\includegraphics{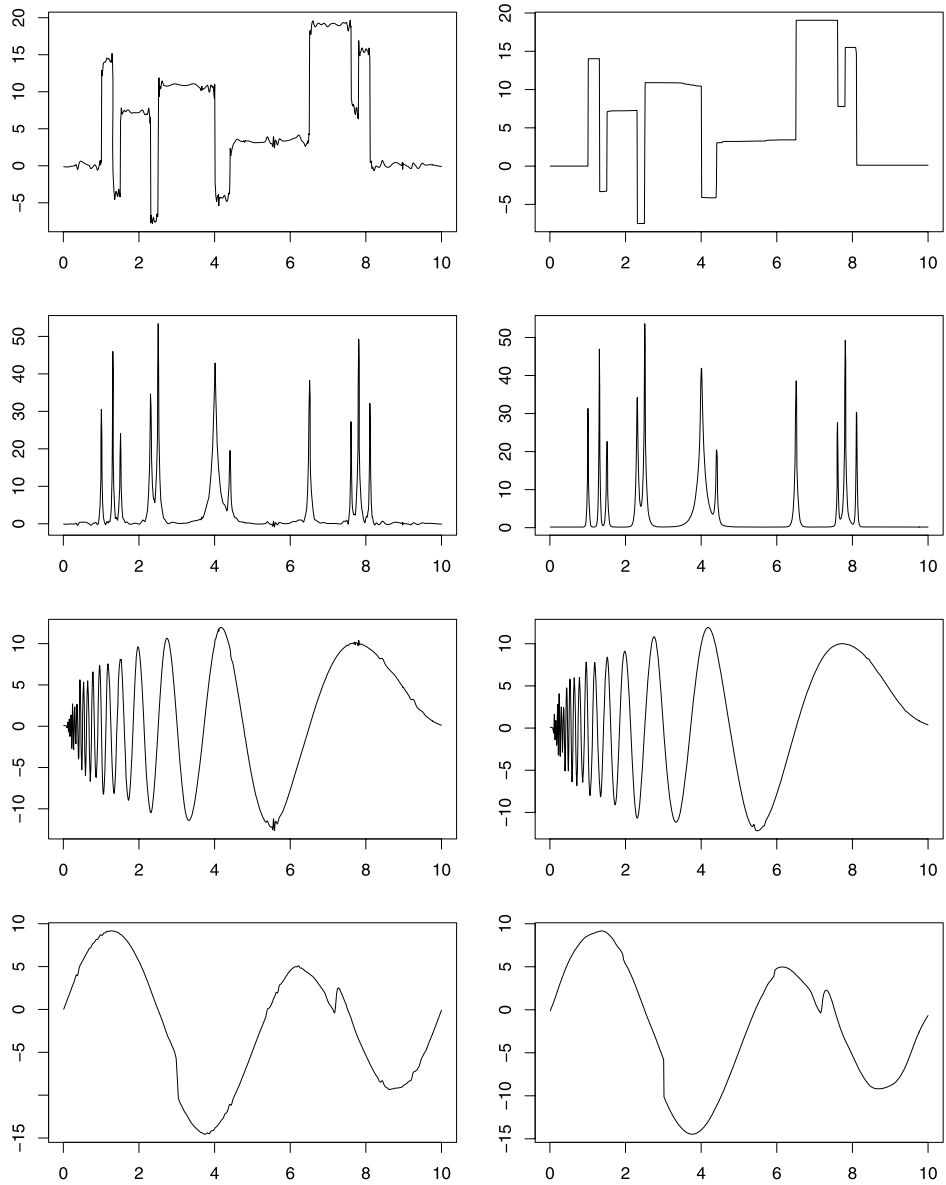}

\caption{Comparison of fitted functions using EBayesThresh
\texttt{beta.laplace} [Johnstone and
Silverman (\protect\citeyear{JohnSilv2005b})] (left column)
and L{\'{e}}vy
Adaptive Regression Kernels (LARK-Gamma) (right column) for the four
test functions. From top to bottom, the test functions are Blocks,
Bumps, Doppler and Heavysine, respectively.}
\label{fsimu}
\vspace*{12pt}
\end{figure}

The performance of each method was measured by its average mean square error
(AMSE), defined as the average value of the MSE given in~(\ref{emse})
over the $100$ replicated simulations. Overall, the performance of the LARK
model is
excellent (Table \ref{tmse}). Both LARK versions generated lower AMSE
values than
did EBayesThresh for all four test functions. LARK also has smaller AMSE
than OCW, except for Doppler, where the methods are comparable. For Blocks,
both LARK and OCW use the Haar wavelet, thus any difference in results
is due to the prior distribution on the function; LARK leads to a 50\%
reduction in AMSE compared to OCW. For the other examples, both OCW and EBayesThresh uses
a Laplace prior distribution for each coefficient in the expansion and the
same wavelet; in all cases it is clear that using a continuous
dictionary is better than the finite-dimensional dictionary (frame) with the nondecimated
wavelets. Lark reconstructions (right column, Figure~\ref{fsimu}) consistently
show less ringing and fewer artifacts than EBayesThresh
(left column).

\section{Applications}\label{sapps}
\subsection{Motorcycle crash data}\label{ssbikes}

\begin{figure}

\includegraphics{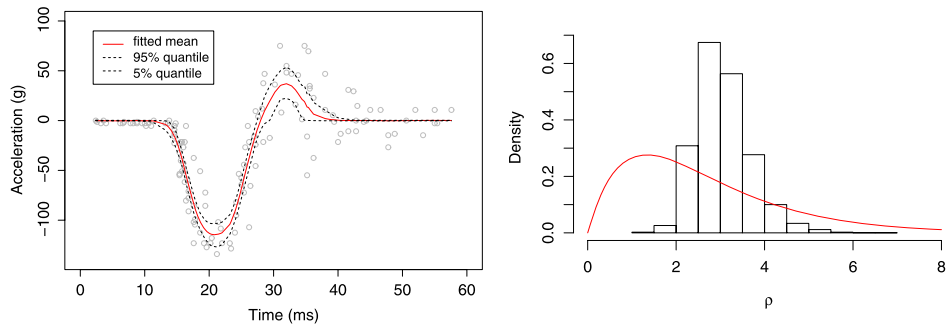}

\caption{Left: results of the LARK model for the motorcycle crash data.
Circles represent the observations; solid line is the posterior mean;
dotted lines are pointwise 90\% Bayesian credible interval for the mean
function. Right: histogram of posterior samples of the exponential power
parameter $\rho$, with prior density (solid line) for
comparison.}\label{fmotor}
\end{figure}

To further illustrate the method, we explore the motorcycle crash experiment
data of \citet{SchmMattSchu1981} considered by \citet{Silv1985},
shown in Figure \ref{fmotor}. The $133$ observations are unequally
spaced, with repeated observations at some time points.
Our focus in this example is to illustrate how a single wide class of
generating functions may be used in LARK, with the data (through the
likelihood) influencing the choice of kernels present in the posterior
distribution.
%
%
We use the power exponential family of kernel functions ${\phi}(x;
\loc, \lambda,
\rho) = \exp\{-\lambda|x-\loc|^\rho\}$, but here (in contrast with the
examples in Section \ref{ssim}) we treat $\rho$ as an uncertain parameter
and make
inference about it from the data. We take the power~$\rho$ to be
common for
all kernels, and use a relatively concentrated Gamma prior distribution
$\rho\sim\Ga(2.0, 0.75)$ with a
50\% HPD interval of $[0.58,2.56]$ which comfortably includes both the
Laplace ($\rho=1$) and Gaussian ($\rho=2)$ kernels as special cases.

The results are summarized in Figure \ref{fmotor}. It is apparent
that the fitted
mean captures the general trend of the data very well, with minimal boundary
effects. The model is parsimonious in the sense we only need $4$
kernels on
average to fit the data. The posterior mean for $\rho$ is
approximately $3$
with most of the posterior mass well above the values ($\rho=1,2$)
for the
Laplace and Gaussian kernels.

%
%
\begin{figure}
\begin{tabular}{@{}cc@{}}

\includegraphics{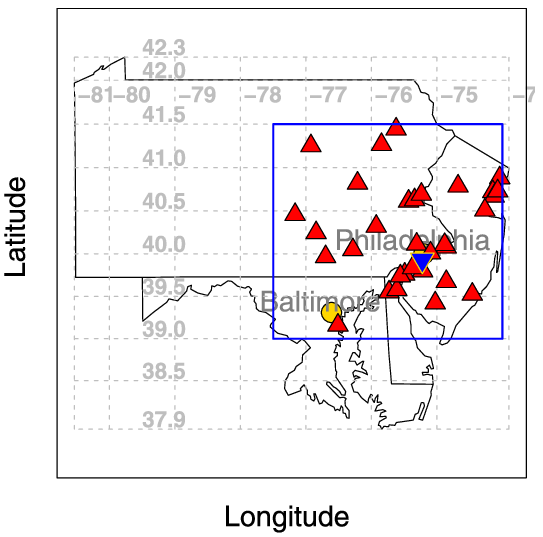}
 & \includegraphics{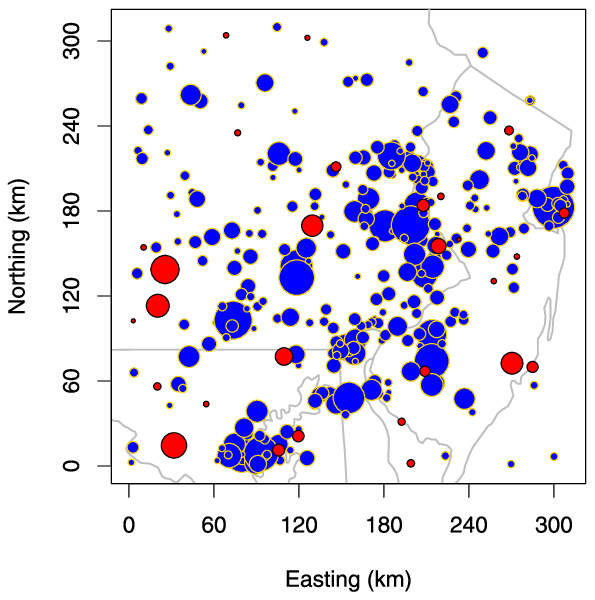}\\
(a) & (b) \\
\end{tabular}
\caption{\textup{(a)} Thirty-three monitors used by EPA to measure hourly
$\mathrm{SO}_2$
concentration in year 2002. The inverted triangle denotes Site 31. The
study area is delineated by a~rectangle that includes parts of
Pennsylvania, Maryland, New Jersey and Delaware, and is blown up in
\textup{(b)} which illustrates locations of kernels from a draw from the posterior
distribution. The blue circles represent aperiodic points and the red
circles represent daily periodic point sources. Circle areas are
proportional to the magnitudes of the point sources they represent.}
\label{fSO2} 
\end{figure}

\subsection{Spatial temporal model}\label{ssair}

In this section, we explore the performance of the LARK approach for modeling
hourly $\mathrm{SO}_2$ concentration levels (measured in ppm) in
Pennsylvania, New
Jersey, Delaware and Maryland [\citet{EPA-AQS2007}]. The
locations of the 33
monitoring stations are shown in Figure \ref{fSO2}; the study region~$\cS$,
delineated by a rectangle in the figure, covers a~$310\km\times310\km
$ area.
We used rescaled coordinates from a Lambert (conformal conic)
projection to
reduce the distortion caused by the earth's curvature. For demonstration
purposes, we restrict analysis to measurements taken during a $144$ hour
period $\cT$ from September of 2002. About $5\%$ of $\mathrm{SO}_2$
readings are
missing (at random) from the data set, which is not a problem for the LARK
model. While Gaussian random field models are popular for modeling
spatial-temporal data, the log transformation typically used in the Gaussian
approach (because the mean function is strictly positive) eliminates
many of
the (important) spiky features of the data. Our Gamma random field prior
distribution allows us to model the data in the original units.

The model can be written in the same simple form as (\ref
{elark-hier}), but
now the $\mathrm{SO}_2$ concentration $Y(x)$ is indexed by points
$x\in\cX
=\cS
\times\cT$ in space--time and the L{\'{e}}vy random measure $\Lmea
(d\omega)$
assigns Gamma-distributed random variables to Borel sets of a space
$\Omega$
of points $\omega=(\sigma,\tau,\Lambda,\lambda)$ that include a location
$(\sigma,\tau)\in\cS\times\cT$ in space--time, a positive-definite
$2\times2$
spatial dispersion matrix $\Lambda\in\cS^+_2$, and a temporal decay rate
$\lambda>0$. We employ a~separable kernel of the form
%
%
\[
{\phi}(x,\omega) = \exp\{-(s-\sigma)'\Lambda(s-\sigma)/2
-\lambda|t-\tau|\}
\]
%
and in the spirit of Higdon [(\citeyear{Higd1998}), Section 3.2]
and Higdon, Swall and Kern [(\citeyear{HigdSwalKern1999}), Section
2.2], we employ a novel parametrization for $\Lambda$ in
terms of its eigenvalues and the orientation of its major axis
[see \citet{Tu2006}, Section~4.2.6, for details on prior
specifications]. In
variations also described in Tu [(\citeyear{Tu2006}), Chapter~4]
accommodation is made
for partial periodicity (due to diurnal patterns associated with daily
variation in ambient temperature, traffic levels, etc.), still within the
framework described by (\ref{elark-hier-all}) but now with more elaborate
choices for $\Omega$ and ${\phi}(x,\omega)$.

The locations of latent point sources from one iteration of the RJ-MCMC
algorithm are presented in Figure \ref{fSO2}(b). Larger latent points
appear to
be clustered in the Baltimore metropolitan area and near the New
Jersey/Penn\-sylvania border. The model's support points are more than a mere
modeling device---they can help analysts identify possible underlying
sources of pollution, or support future decisions on monitor locations.

\begin{figure}[b]

\includegraphics{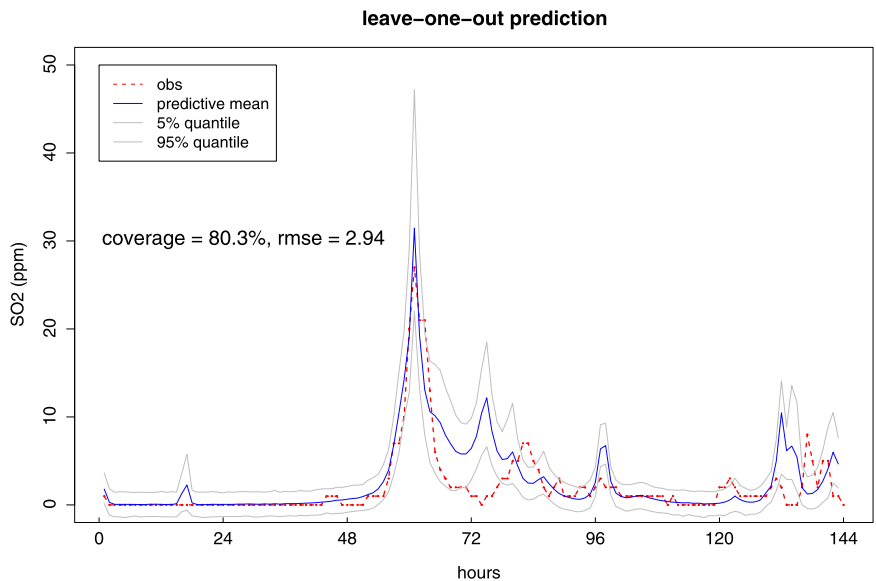}

\caption{Out-of-sample predictions for Site $31$. Dashed line represents
observed time series, solid line represents predictive mean curve.
Gray lines are $90\%$ posterior predictive intervals.}
\label{fSO2-pred}
\end{figure}

The predictive power of the model is validated through out-of-sample
prediction. The model was fit excluding data from Site $31$ [the inverted
triangle in Figure \ref{fSO2}(a)], and then its predictions were
compared with
reported measurements from that site for the entire $144$ hours. The result
shown in Figure \ref{fSO2-pred} is promising. The major peak was captured
%
%
clearly, and $90\%$ pointwise Bayesian credible intervals cover in
excess of
$80\%$ of the true observations. This was a challenging out-of-sample
prediction problem due to low cross-correlations among sites. We are
currently refining features of the prior distributions to incorporate known
point sources.

\section{Discussion}\label{sdisc}
In this article, we have developed a fully Bayesian adaptive kernel method,
LARK, for nonparametric function estimation. The LARK model is based on a
stochastic expansion of functions in a continuous overcomplete dictionary,
and may be expressed as a stochastic integral of a kernel or other generating
function with respect to a L{\'{e}}vy random field. When (\ref
{el1bound}) is
satisfied (so compensation is unnecessary), the L{\'{e}}vy field is a~random
signed measure. By using a \textit{positive} random measure and positive
kernel family, LARK models provide natural constructions for nonnegative
functions (as in Section \ref{ssair}); with \textit{signed} measures,
unconstrained
functions may be modeled (as in Sections \ref{ssim} and \ref
{ssbikes}). The kernel
parameters are location-specific and thus adapt to local features of the
data. As with wavelets, the adaptive smoothing using LARK preserves local
features such as discontinuities and high peaks and is especially
useful for
modeling inhomogeneous functions. The LARK approach does not require that
the data be equally-spaced without missing observations nor that the sample
size be a dyadic power as is a commonly required of many wavelet
methods.\looseness=1

The RJ-MCMC algorithm developed for fitting LARK provides an\break automatic
stochastic search mechanism for finding sparse representations of a~function.
The algorithm is computationally efficient [requiring only
\mbox{$O(n\cdot M)$}~operations for data including $n$ observations and an MCMC stream of length~$M$],
as dictionary elements are calculated only when needed. Kernel methods
such as Support Vector Machines (SVMs) and Bayesian Relevance Vector Machines
[or RVMs, \citet{Tipp2001}] employ all data points as kernel
locations, but
attain sparsity by shrinking coefficients to zero. LARK provides additional
flexibility by not restricting kernel locations. Many competing sparse
methods, including the Dantzig Selector and Lasso, require the
a~priori selection of a pre-specified number of dictionary elements.
Evaluating these kernels on a sufficiently fine grid will exceed the
computational cost of LARK. Fine grids also lead to extreme
multicollinearity in these approaches, that may lead both to numerical
instability and violation of the conditions needed for sparse solutions.



\vspace*{6pt}\subsection{Extensions}\label{ssextend}

It is straightforward to implement LARK with wide classes of generating
functions including wavelets, structural elements in texture analysis, and
splines. Unlike support vector machines or other methods based on Mercer
kernels [\citet{PillWuetal2007}], the LARK approach does not require
symmetry, continuity or simple functional forms. While it is often
convenient to use kernels based on some distance metric, arbitrary generating
functions may be tailored to the problem at hand as illustrated in the
space--time example of Section \ref{ssair}. The LARK modeling approach adapts
readily to problems in any number of dimensions.

In Section \ref{sspaces}, we present conditions for LARK models to
belong to
the same
Besov space as their generating functions, for L{\'{e}}vy measures and
generating
functions that satisfy the stringent local $L_1$-bound of (\ref{enocomp}).
In the more general case, where (\ref{enocomp}) fails and
compensation is
required, we are able to establish similar results only for $\Besov$
with $p
= q = 2$ (equivalent to $\Sob$).
We are exploring extensions to the general case, but the additional drift
term that arises in compensation complicates confirming the convergence
of~$f_\eps$ to $f$ in $\Besov$ for general $p,q$.


Work is also on-going in establishing conditions for posterior consistency
for function estimation. Extending methods of \citet{ChouGhosRoy2004},
\citet{GhosVaar2007a} and \citet{ChoiSche2007},
\citet{Pill2008} has
verified posterior consistency for certain LARK models with Gaussian
measurement errors in work that will be reported elsewhere.

\begin{appendix}\label{app}
\vspace*{3pt}\section{Details of proofs}\vspace*{3pt}
\label{aproofs}

\begin{prop}\label{propdiff}
For a function $\G(\cdot)\in L_p(\bbR^d)$ and its scaled translate
$\G(\Scale( \cdot- \loc))$ with $\loc\in\bbR^d$ and positive
definite matrix $\Scale\in\LamS$, the $L_p$ norm of $\G(\Scale
( \cdot-
\loc))$ and the $L_p$ norm of its $m$th forward differences are given
by
%
\begin{equation} \label{eqdiffnorm}
\bigl\| \G\bigl(\Scale( \cdot- \loc)\bigr)\bigr\|_p = |\Scale|^{1/p} \| \G\|_p
\bigl\| \Delta^m_h \G\bigl(\Scale( \cdot- \loc)\bigr)\bigr\|_p = |\Scale
|^{1/p} \| \Delta^m_{\lambda h} \G\|_p,\hspace*{-28pt}
\end{equation}
where $|\Scale|$ denotes the determinant of $\Scale$.
\end{prop}
\begin{pf}
By a change of variables $\loc\mapsto u=\Scale(x-\loc)$,
\begin{eqnarray*}
\bigl\| \Delta^{m}_h \G\bigl(\Scale( \cdot- \loc)\bigr)\bigr\|_p
& = & \biggl\{\int\bigl|\Delta
^{m}_h \G\bigl(\Scale(x - \loc) \bigr)\bigr|^p \,dx \biggr\} ^{1/p} \\
& = & \Biggl\{\int\Biggl| \sum_{k=0}^m \pmatrix{m\cr k} (-1)^{m - k} \G\bigl(\Scale(
x + k h - \loc)\bigr) \Biggr|^p \,dx \Biggr\}^{1/p} \\
& = &|\Scale|^{-1/p} \Biggl\{\int\Biggl| \sum_{k=0}^m \pmatrix{m\cr k} (-1)^{m -
k} \G(u+k\Scale h) \Biggr|^p \,du \Biggr\}^{1/p} \\
& = &|\Scale|^{-1/p} \| \Delta^m_{\Scale h} \G\|_p.
\end{eqnarray*}
The proof for the $L_p$ norm of $\G(\Scale( \cdot- \loc))$
follows by the same change of variables.
\end{pf}

\subsection{\texorpdfstring{Proof of Lemma \protect\ref{lASS}}{Proof of Lemma 1}}
\label{alem1}
First, consider the case $b > 1$ and $a \in\bbR$. Then
\begin{eqnarray*}
&&\iiint_{\bbR\times[1, \infty)\times\bbT}
\bigl(1 \wedge|\Z\GB(u)^r | \lambda^{-a}\bigr) \lambda^{-b}
\pi_\Z(d\Z) \,d\lambda \,du \\
&&\qquad <
\iiint_{\bbR\times[1, \infty)\times\bbT}
\lambda^{-b}
\pi_\Z(d\Z) \,d\lambda \,du \\
&&\qquad = \frac{1}{b-1}<\infty.
\end{eqnarray*}
Next, consider the case of $b<1$ and $a> 1- b$ (which imply $a>0$):
\begin{eqnarray*}
&&\iiint_\RRT
\bigl(1 \wedge|\Z\GB(u)^r| \lambda^{-a}\bigr) \lambda^{-b}
\pi_\Z(d\Z) \,d\lambda \,du\\[2pt]
&&\qquad= \iint_{|z\GB(u)^r|>1} 
\int_{1}^{|\Z\GB(u)^r |^{1/a}} 
\lambda^{-b} \,d \lambda\pi_\Z(d\Z) \,du \\[2pt]
&&\qquad\quad{}+ \iint_{\bbR\times\bbT} |\Z\GB(u)^r|
\int_{1\vee|\Z\GB(u)^r |^{1/a}}^{\infty} \lambda^{-a-b} \,d \lambda
\pi_\Z(d\Z) \,du \\[2pt]
&&\qquad= \iint_{|z\GB(u)^r|>1} 
\frac{\lambda^{1 - b}}{1 - b}
\bigg|_{\lambda=1}^{\lambda=|\Z\GB(u)^r |^{1/a}} \pi_\Z(d\Z)
\,du \\[2pt]
&&\qquad\quad{}+
\iint_{\bbR\times\bbT} |\Z\GB(u)^r| \frac{\lambda^{(1 - a - b)}}
{1 - a - b} \bigg|_{\lambda=1\vee|\Z\GB(u)^r |^{1/a}}
^{\lambda=\infty} \pi_\Z(d\Z) \,du \\[2pt]
&&\qquad= \iint_{|\Z\GB(u)^r|>1}
\frac{1-| \Z\GB(u)^r|^{({1 - b})/{a}}}
{b - 1} \pi_\Z(d\Z) \,du\\[2pt]
&&\qquad\quad{}+ \iint_{|\Z\GB(u)^r|>1}\frac{| \Z\GB(u)^r|^{({1 - b})/{a}}}
{a + b -1} \pi_\Z(d\Z) \,du \\[2pt]
&&\qquad\quad{} + \iint_{|\Z\GB(u)^r|\le1}
\frac{|\Z\GB(u)^r|} {a + b -1} \pi_\Z(d\Z) \,du \\[2pt]
&&\qquad\le\frac{1}{b - 1} +
\iint_{\bbR\times\bbT} | \Z\GB(u)^r|^{({1 - b})/{a}}
\frac{ a}{(a+b-1)(1 - b)} \pi_\Z(d\Z) \,du \\[2pt]
&&\qquad= \frac{1}{b - 1} +
\frac{a} {(a + b -1)(1 - b)}
\int_{\bbR} |\Z|^{({1 - b})/{a}} \pi_\Z(d\Z)
\int_{\bbT} |\GB(u)^r|^{({1 - b})/{a}} \,du
\\[2pt]
&&\qquad<\infty
\end{eqnarray*}
for $a+b>1$ if $r=0$, and for $ap+b\ge1$ if $r=1$ [since
$\GB\in L^*_p(\bbT)$], which is implied by $a>1-b$.

Now consider the case of $b=1$ and $a>0$:
\begin{eqnarray*}
&&\iiint_\RRT
\bigl(1 \wedge|\Z\GB(u)^r| \lambda^{-a}\bigr) \lambda^{-b}
\pi_\Z(d\Z) \,d\lambda \,du\\
&&\qquad= \iint_{|z\GB(u)^r|>1} 
\int_{1}^{|\Z\GB(u)^r |^{1/a}} 
\lambda^{-1} \,d \lambda\pi_\Z(d\Z) \,du \\[2pt]
&&\qquad\quad{}+ \iint_{\bbR\times\bbT} |\Z\GB(u)^r|
\int_{1\vee|\Z\GB(u)^r |^{1/a}}^{\infty} \lambda^{-a-1} \,d \lambda
\pi_\Z(d\Z) \,du \\[2pt]
&&\qquad= \iint_{|z\GB(u)^r|>1} 
\log\lambda
\bigg|_{\lambda=1}^{\lambda=|\Z\GB(u)^r |^{1/a}} \pi_\Z(d\Z)
\,du \\[2pt]
&&\qquad\quad{}+ \iint_{\bbR\times\bbT} |\Z\GB(u)^r| \frac{\lambda^{-a}}
{- a } \bigg|_{\lambda=1\vee|\Z\GB(u)^r |^{1/a}}
^{\lambda=\infty} \pi_\Z(d\Z) \,du \\[2pt]
&&\qquad= \iint_{|\Z\GB(u)^r|>1}
{\frac{1}{a}\log}| \Z\GB(u)^r| \pi_\Z(d\Z) \,du
+ \iint_{|\Z\GB(u)^r|>1}\frac1a \pi_\Z(d\Z) \,du \\[2pt]
&&\qquad\quad{} + \iint_{|\Z\GB(u)^r|\le1}
\frac{|\Z\GB(u)^r|} {a} \pi_\Z(d\Z) \,du \\[2pt]
&&\qquad\le\frac1a\iint_{\bbR\times\bbT} {\log_+}| \Z\GB(u)^r|
\pi_\Z(d\Z) \,du + \frac1a \\[2pt]
&&\qquad<\infty
\end{eqnarray*}
since $\log_+(z\GB^r)=(0\vee\log|z\GB^r|) \le|z|+|\GB|^r$ and
$\GB\in L^*_1(\bbT)$.

\vspace*{3pt}\subsection{\texorpdfstring{Proof of Theorem \protect\ref{tconverge}}{Proof of Theorem 2}}
\label{athm2}

For any compensator function $\cmp$ satisfying~(\ref{ecomp}) there are numbers $c_j\in(0,\infty)$ such that
\[
|\cmp|\le c_0,\qquad
|\beta-\cmp|\le c_1 (|\beta|\wedge\beta^2),\qquad
|\cmp| \le c_2 ( 1 \wedge|\beta|)
\]
for all $\beta\in\bbR$. Fix $0<\eps\le1$ and a function $\phi
\dvtx\bbR\times
\Omega\to\bbR$ satisfying (\ref{eexist}); let~$B_a$, $B_b$ and $B_c$
be the
values of the integrals from (\ref{eexa})--(\ref{eexc}),
respectively.
%
To complete the proof of Theorem \ref{tconverge} it suffices to show that
each of
the two terms from (\ref{epoi-comp-res}),
%
\begin{eqnarray}\label{ebnds}
X&\equiv&\iint_{N_\eps}\bigl(\beta- \cmp\bigr)
\phi(\omega) \Np(\dbdo)
\quad\mbox{and}\nonumber\\[-8pt]\\[-8pt]
Y &\equiv&\iint_{N_\eps} \cmp\phi(\omega)
\Nt(\dbdo),\nonumber
\end{eqnarray}
converges to zero in probability as $\eps\to0$.
Write the first integral in (\ref{ebnds}) as the sum of two parts:
\[
X \equiv \iint_{N_\eps} \bigl(\beta- \cmp\bigr)
\phi(\omega) \Np(\dbdo)
= X_1 + X_2
\]
with
\begin{eqnarray*}
X_1 &\equiv&\iint_{{ N_\eps\cap
[ |\beta\phi|\le1 ]}}
\bigl(\beta- \cmp\bigr) \phi(\omega) \Np(\dbdo),\\[2pt]
X_2 &\equiv&\iint_{{ N_\eps\cap
[|\beta\phi|> 1 ]}}
\bigl(\beta- \cmp\bigr) \phi(\omega) \Np(\dbdo).
\end{eqnarray*}
Then
\begin{eqnarray*}
\E|X_1| &\le& c_1 \iint_{{ N_\eps\cap[ |\beta\phi|\le1 ]}}
(|\beta|\wedge\beta^2) |\phi(\omega)| \nu(\dbdo)
\\[2pt]
&=& c_1 \iint_{{ N_\eps\cap
[ |\beta|\le1] \cap[ |\beta\phi|\le1 ]}}
(1\wedge\beta^2) |\phi(\omega)| \nu(\dbdo)\\[2pt]
&&{}+ c_1 \iint_{{ N_\eps\cap
[ |\beta|>1] \cap[ |\beta\phi|\le1 ]}}
\bigl(1\wedge|\beta\phi(\omega)|\bigr) \nu(\dbdo)\\[2pt]
&\le& c_1 (B_c + B_a) < \infty,
\end{eqnarray*}
so $X_1\to0$ in $L_1$ as $\eps\to0$ by Lebesgue's dominated convergence
theorem since the indicator function $\mathbf{1}_{\{N_\eps\}}(\beta
,\omega)$
tends to
zero a.e. $(\nu)$ as $\eps\to0$. Now consider~$X_2$:
\begin{eqnarray*}
\nu\bigl(\{(\beta,\omega)\dvtx|\beta\phi(\omega)|>1\}\bigr)
&=& \iint_{{
[ |\beta|\le1] \cap[ |\beta\phi|> 1 ]}}
1 \nu(\dbdo)\\[2pt]
&&{}
+ \iint_{{ [ |\beta|>1] \cap[ |\beta\phi|> 1 ] }}
1 \nu(\dbdo)\\[2pt]
&\le&\iint_{{
[ |\beta|\le1] \cap[ |\beta\phi|> 1 ]}}
\bigl(|\beta\phi(\omega)|\wedge|
\beta\phi(\omega)|^2\bigr)
\nu(\dbdo) \\[2pt]
&&{}
+ \iint_{{ [ |\beta|>1] \cap[ |\beta\phi|> 1 ] }}
\bigl(1\wedge|\beta\phi(\omega)|\bigr) \nu(\dbdo)\\[2pt]
&\le& B_b + B_a < \infty,
\end{eqnarray*}
so almost surely the random support of $\Np(\dbdo)$ in
$[|\beta\phi|>1]$ is a finite set disjoint from $\bigcap_{\eps>0}
N_\eps$; it
follows that $\Np(N_\eps\cap[ |\beta\phi(\omega)|>1 ]
)\to0$
and hence $X_2\to0$ almost surely as $\eps\to0$.

Similarly, we write the second integral in (\ref{ebnds}) as the sum
of four
parts:
\[
Y \equiv\iint_{N_\eps} \cmp\phi(\omega) \Nt(\dbdo)
= Y_1 + Y_2 + Y_3 + Y_4
\]
with
\begin{eqnarray*}
Y_1 &\equiv& \iint_{{ N_\eps\cap
[ |\beta|\le1]\cap
[ |\beta\phi|\le1 ] }}\cmp\phi(\omega) \Nt(\dbdo),\\
Y_2 &\equiv&\iint_{{ N_\eps\cap
[ |\beta|\le1]\cap
[ |\beta\phi|> 1 ] }}\cmp\phi(\omega) \Nt(\dbdo),\\
Y_3 &\equiv&\iint_{{ N_\eps\cap
[ |\beta|>1]\cap
[ |\beta\phi|\le1 ] }}\cmp\phi(\omega) \Nt(\dbdo),\\
Y_4 &\equiv&\iint_{{ N_\eps\cap
[ |\beta|>1]\cap
[ |\beta\phi|> 1 ] }}\cmp\phi(\omega) \Nt(\dbdo).
\end{eqnarray*}
Now
\begin{eqnarray*}
\E|Y_1|^2 &=& \iint_{{ N_\eps\cap
[ |\beta|\le1]\cap
[ |\beta\phi|\le1 ] }}
\cmp^2 \phi(\omega)^2 \nu(\dbdo)\\
&\le& c_2^2
\iint_{{ N_\eps\cap
[ |\beta|\le1]\cap
[ |\beta\phi|\le1 ] }}
|\beta\phi(\omega)|^2 \nu(\dbdo)\\
&=& c_2^2
\iint_{{ N_\eps\cap
[ |\beta|\le1]\cap
[ |\beta\phi|\le1 ] }}
\bigl(|\beta\phi(\omega)|\wedge
|\beta\phi(\omega)|^2 \bigr) \nu(\dbdo)\\
&\le& c_2^2 B_b<\infty,
\end{eqnarray*}
so $Y_1\to0$ in $L_2$ (and hence also in $L_1$) as $\eps\to0$ by LDCT,
\begin{eqnarray*}
Y_2 &\equiv& \iint_{{ N_\eps\cap
[ |\beta|\le1]\cap
[ |\beta\phi|> 1 ] }}
\cmp\phi(\omega) \Nt(\dbdo)\\
&=& \iint_{{ N_\eps\cap
[ |\beta|\le1]\cap
[ |\beta\phi|> 1 ] }}
\cmp\phi(\omega) \Np(\dbdo) \\
&&{}-
\iint_{{ N_\eps\cap
[ |\beta|\le1]\cap
[ |\beta\phi|> 1 ] }}
\cmp\phi(\omega) \nu(\dbdo),\\
\E|Y_2| &\le& 2 \iint_{{ N_\eps\cap
[ |\beta|\le1]\cap
[ |\beta\phi|> 1 ] }}
|\cmp| |\phi(\omega)| \nu(\dbdo)\\
&\le& 2 c_2 \iint_{{ N_\eps\cap
[ |\beta|\le1]\cap
[ |\beta\phi|> 1 ] }}
|\beta\phi(\omega)| \nu(\dbdo)\\
&=& 2 c_2 \iint_{{ N_\eps\cap
[ |\beta|\le1]\cap
[ |\beta\phi|> 1 ] }}
\bigl(|\beta\phi(\omega)|\wedge
|\beta\phi(\omega)|^2 \bigr) \nu(\dbdo)\\
&\le& 2 c_2 B_b<\infty,
\end{eqnarray*}
so $Y_2\to0$ in $L_1$ as $\eps\to0$ by dominated convergence,
\begin{eqnarray*}
Y_3 &\equiv&\iint_{{ N_\eps\cap
[ |\beta|>1]\cap
[ |\beta\phi|\le1 ] }}
\cmp\phi(\omega) \Nt(\dbdo)\\[2pt]
&=& \iint_{{ N_\eps\cap
[ |\beta|>1]\cap
[ |\beta\phi|\le1 ] }}
\cmp\phi(\omega) \Np(\dbdo) \\[2pt]
&&{} -
\iint_{{ N_\eps\cap
[ |\beta|>1]\cap
[ |\beta\phi|\le1 ] }}
\cmp\phi(\omega) \nu(\dbdo),\\[2pt]
\E|Y_3| &\le& 2 \iint_{{ N_\eps\cap
[ |\beta|>1]\cap
[ |\beta\phi|\le1 ] }}
|\cmp| |\phi(\omega)| \nu(\dbdo)\\[2pt]
&\le& 2 c_2 \iint_{{ N_\eps\cap
[ |\beta|>1]\cap
[ |\beta\phi|\le1 ] }}
|\beta\phi(\omega)| \nu(\dbdo)\\[2pt]
&=& 2 c_2 \iint_{{ N_\eps\cap
[ |\beta|>1]\cap
[ |\beta\phi|\le> 1 ] }}
\bigl(1 \wedge|\beta\phi(\omega)|
\bigr) \nu(\dbdo)\\[1pt]
&\le& 2 c_2 B_a<\infty,
\end{eqnarray*}
so $Y_3\to0$ in $L_1$ as $\eps\to0$. Finally, for $Y_4$,
\begin{eqnarray*}
Y_4 &\equiv& \iint_{{ N_\eps\cap
[ |\beta|>1]\cap
[ |\beta\phi|> 1 ] }}
\cmp\phi(\omega) \Nt(\dbdo)\\[2pt]
&=& \iint_{{ N_\eps\cap
[ |\beta|>1]\cap
[ |\beta\phi|> 1 ] }}
\cmp\phi(\omega) \Np(\dbdo)\\[2pt]
&&{}-
\iint_{{ N_\eps\cap
[ |\beta|>1]\cap
[ |\beta\phi|> 1 ] }}
\cmp\phi(\omega) \nu(\dbdo),\\[2pt]
\E|Y_4| &\le& 2
\iint_{{ N_\eps\cap
[ |\beta|>1]\cap
[ |\beta\phi|> 1 ] }}
|\cmp\phi(\omega)| \nu(\dbdo)\\[2pt]
&\le& 2 c_0 \iint_{{ N_\eps\cap
[ |\beta|>1]\cap
[ |\beta\phi|> 1 ] }}
|\phi(\omega)| \nu(\dbdo)\\[2pt]
&\le& 2 c_0 \iint_{{ N_\eps\cap
[ |\beta|>1]\cap
[ |\beta\phi|> 1 ] }}
(1\wedge\beta^2) |\phi(\omega)| \nu(\dbdo)\\[2pt]
&\le& 2 c_0 B_c<\infty
\end{eqnarray*}
so $Y_4\to0$ in $L_1$ as $\eps\to0$, completing the proof of
Theorem \ref{tconverge}.

\vspace*{2pt}\section{Reversible-jump MCMC procedures}\vspace*{2pt}\label{arjmcmc}

A typical RJ-MCMC procedure for sampling varying-dimensional parameters
involves at least three types of moves (Birth, Death and Update); we use
Metropolis--Hastings steps for each of these. Our trans-dimensional update
steps entail altering the value $(\beta_j^*, \omega_j^*)$ of one point
$(\beta_j,\omega_j)$. We select\vadjust{\goodbreak} $j\sim\Un(0\dvtx J-1)$ for proposed updating,
then take Gaussian random walk steps successively in the coefficient
$\beta_j$, the location parameter $\loc_j$, and the log kernel shape
parameter, $\log\lambda_j$. Step sizes are chosen to achieve approximately
$30\%$ acceptance rates for each class of updates. One novel feature is that
when the proposed update of some coefficient $\beta_j$ falls\vspace*{1pt} in the truncated
region $\beta^*_j\eta\in(-\eps,\eps)$, the move is treated as a
Death, the
point $(\beta_j ,\omega_j)$ is removed and $J$ is decremented. This is
advantageous as it automatically focuses on small magnitude
coefficients for
removal (rather than a random selection as in the typical RJ-MCMC Death
step). A Birth step entails generating a new point $(\beta^*,\omega
^*)$ to
be included among the $\{(\beta_j,\omega_j)\}$ and incrementing $J$
by one.
We use a double exponential birth distribution with rate~$\eta/\eps$,
conditioned to exceed $|\beta_j|\eta>\eps$ so that proposed
coefficients are
small, balancing the ``Death'' of small coefficients in the Update step to
attain the target acceptance rates. The fixed-dimensional parameters are
sampled using a conventional Metropolis--Hastings approach
[\citet{GilkRichSpie1996}, Section 1.3.3].
Each of these inexpensive update steps requires only $O(n)$ operations [in
contrast to Gaussian methods, which may require $O(n^3)$], so the method
scales well in the number $n$ of observations.
Further details of the RJ-MCMC are available in [\citet{Tu2006},
Appendix~A.1, pages 116 and 117]. An \texttt{R} package [R~Development Core Team
(\citeyear{R2004})] implementing LARK
is under development by the authors and will be made publicly
available.

\section{Examples of LARK prior realizations}\label{aegs}
%

\begin{figure}[h]
\begin{tabular}{@{}c@{\hspace*{1pt}}c@{}}

\includegraphics{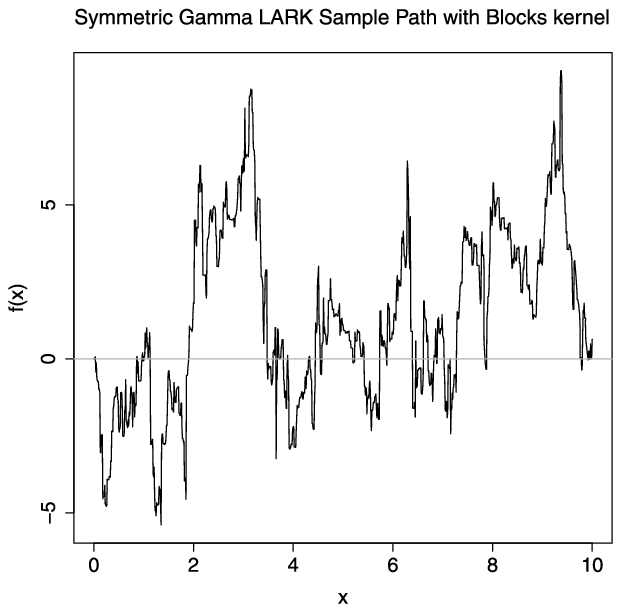}
 & \includegraphics{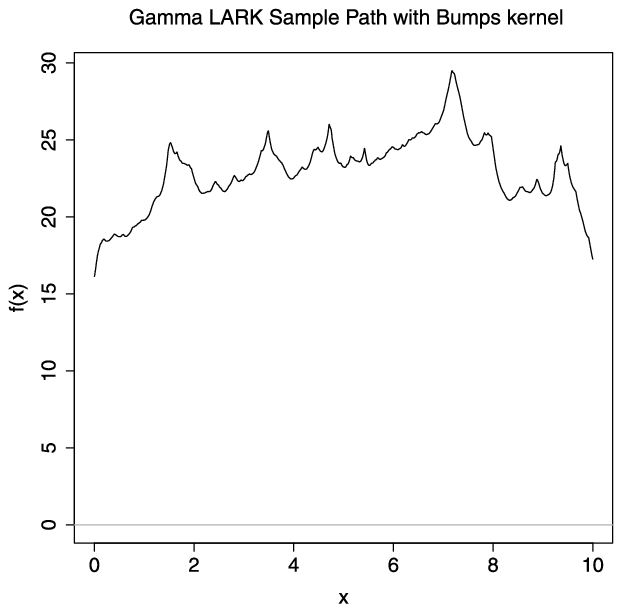} \\
(a)&(b)
\end{tabular}
\caption{Four realizations from LARK prior distribution with \textup{(a)} Blocks
kernel and Symmetric Gamma L{\'{e}}vy measure; \textup{(b)} Bumps kernel and Gamma
L{\'{e}}vy measure; \textup{(c)}, \textup{(d)} Doppler kernel and Cauchy L{\'{e}}vy measure, with
$J=1000$ for \textup{(a)--(c)} and $J=10$ for \textup{(d)} components.
Hyperparameters $a_\lambda$, $b_\lambda$, $a_\gamma$, $b_\gamma$, $a_\eta$,
$b_\eta$ and
$\eps$ are given in Table \protect\ref{thyp}.}\label{fRealize}
\end{figure}

\setcounter{figure}{4}
\begin{figure}
\begin{tabular}{@{}c@{\hspace*{1pt}}c@{}}

\includegraphics{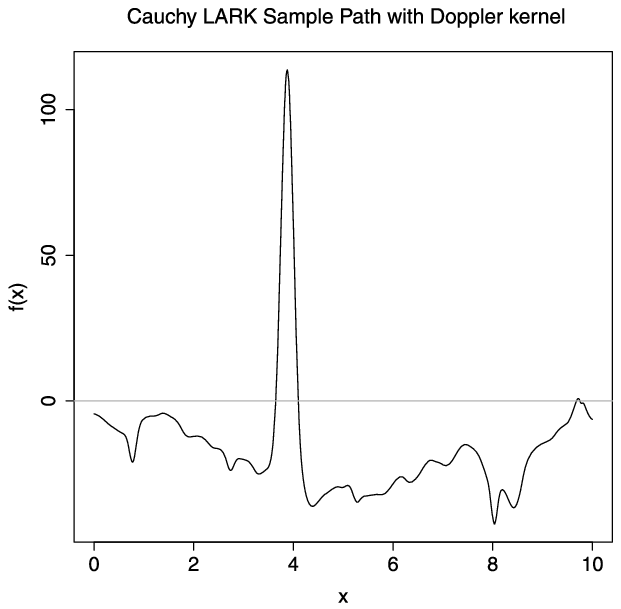}
 & \includegraphics{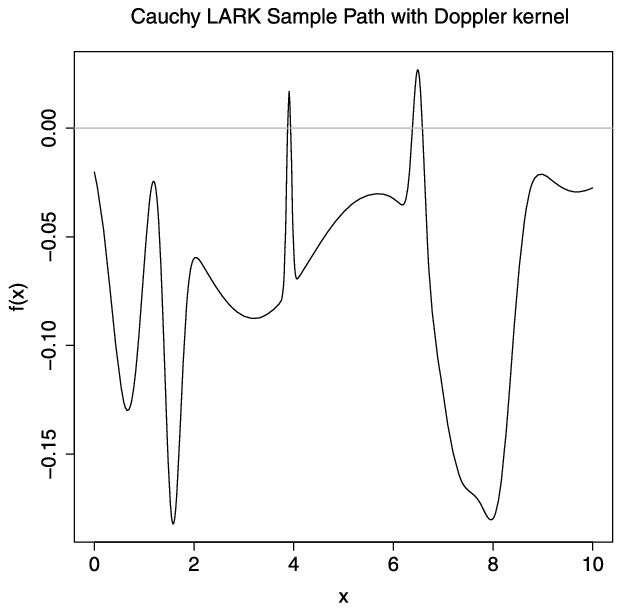}\\
(c)&(d)
\end{tabular}
\caption{(Continued.)}
\end{figure}
\mbox{}\vadjust{\goodbreak}
\end{appendix}

\section*{Acknowledgments}
The authors would like to thank Natesh Pillai, three
referees, the Associate Editor and the Editor for helpful comments and
suggestions.


%

%
\printaddresses

\end{document}